\documentclass[12pt,reqno]{amsart}
\usepackage{amsmath}%
\usepackage{amsfonts}%
\usepackage{amssymb}%
\usepackage[margin=.5in]{geometry}
\usepackage{mathtools}
\usepackage{graphicx, caption, subcaption}
\usepackage{array,multirow}
\usepackage{ifpdf}
\ifpdf
\usepackage{graphicx}
\else
\usepackage[draft]{graphicx}
\fi
\usepackage{color}
\usepackage{enumerate}
\usepackage[T1]{fontenc}
\usepackage{booktabs,multirow}
\usepackage{enumerate}
\usepackage{listings}
\usepackage{color}
\usepackage{caption}
\usepackage{subcaption}
\usepackage{xcolor}
\usepackage{bm}
\usepackage[hyphens]{url}
\usepackage{hyperref} 

\hypersetup{
    colorlinks=true,
    linkcolor=blue,
    filecolor=magenta,      
    urlcolor=cyan,
}
\definecolor{dkgreen}{rgb}{0,0.6,0}
\definecolor{gray}{rgb}{0.5,0.5,0.5}
\definecolor{mauve}{rgb}{0.58,0,0.82}

\newcommand{\Arg}{\text{Arg}}
\newcommand{\pLog}{\texttt{pLog}}
\newcommand{\branchF}{\text{branchF}}
\newcommand{\leafF}{\text{leafF}}
\newcommand{\Log}{\text{ Log}}
\newcommand{\quotes}[1]{``#1''}
\begin{document}

\title[]{A systematic approach to computing and indexing fixed points of an  iterated exponential}
\author{Dominic C. Milioto\\ \MakeLowercase{icorone@hotmail.com}}
\email{icorone@hotmail.com}
\date{\today}
%
\begin{abstract}
This paper describes a systematic method of numerically computing and indexing fixed points of $z^{z^w}$ for fixed $z$ or equivalently, the roots of $T_2(w;z)=w-z^{z^w}$.  The roots are computed using a modified version of fixed-point iteration and indexed by integer triplets $\{n,m,p\}$ which associate a root to a unique branch of $T_2$. This naming convention is proposed sufficient to enumerate all roots of the function with $(n,m)$   enumerated by $\mathbb{Z}^2$.  However, branches near the origin can have multiple roots.  These cases are identified by the third parameter $p$.  This work was done with rational or symbolic values of $z$ enabling arbitrary precision arithmetic.  A selection of roots up to order $\{10^{12},10^{12},p\}$ with $|z|\leq 10^{12}$ was used as test cases.  Results were accurate to the precision used in the computations, generally between $30$ and $100$ digits. Mathematica ver. $12$ was used to implement the algorithms.  
\end{abstract}
\maketitle
\pagestyle{plain}
\section{Introduction}
A method is developed for computing and indexing fixed points of $z^{z^w}$ for fixed $z$ or equivalently, the zeros of the function  $T_2(w;z)=w-z^{z^w}$ for complex $z$ and $w$  with $z=re^{i\alpha}$.  Values of $r$ and $\alpha$ were restricted to rational or symbolic quantities allowing for arbitrary precision arithmetic. 

For constant $z$, define
\begin{equation}
\begin{aligned}
T_1(w;z)&=w-z^w\equiv w-e^{w\Log(z)}, \\
T_2(w;z)&=w-z^{z^w}\equiv w-e^{z^w \Log(z)}\equiv e^{\Log(z)e^{w\Log(z)}},
\end{aligned}
\label{equation:eqn0}
\end{equation} 
where $\Log(z)=\ln|z|+i\theta$ with $\theta=\Arg(z)$ and $-\pi<\Arg(z)\leq \pi$.  $T_1$ is called the $1$-cycle expression and $T_2$, the $2$-cycle expression.  An explicit formula  for the zeros of $T_1$ exits.  However, in this paper $T_1$ is analyzed without it so that the analysis of $T_2$ easily follows.   

\begin{figure}[!ht]
     \centering
     \begin{subfigure}{0.49\textwidth}
         \centering
         \includegraphics[width=0.75\textwidth]{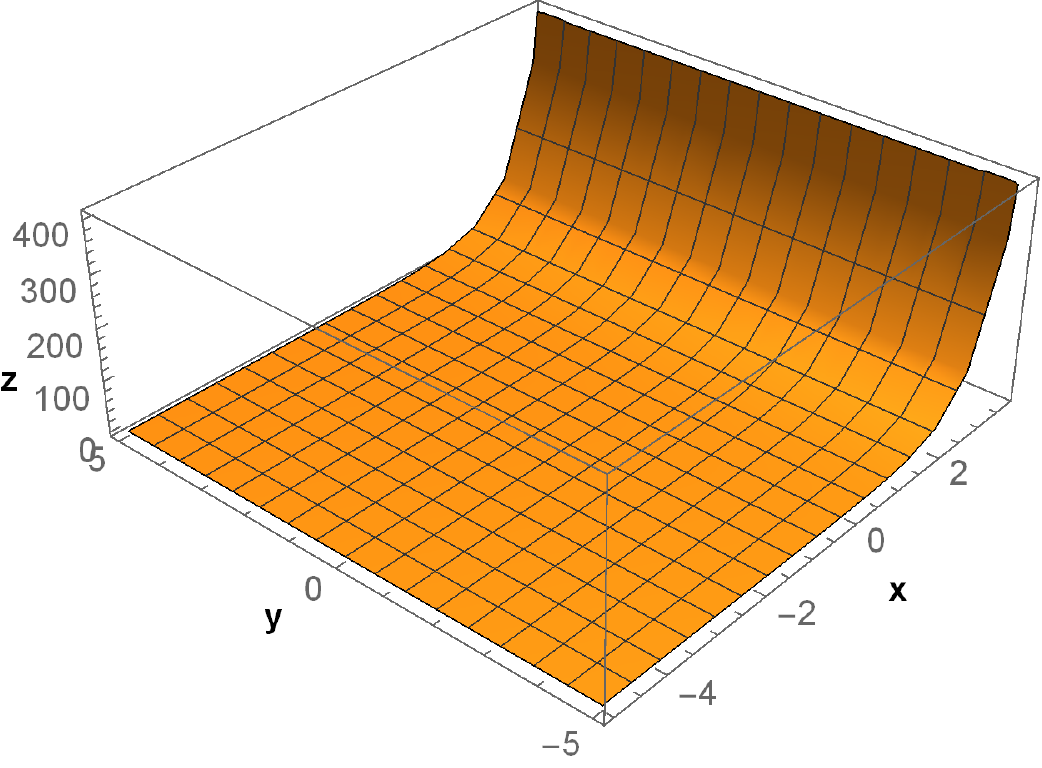}
         \caption{$e^{ax}$}
         \label{figure:figure0a}
     \end{subfigure}
     \hfill
     \begin{subfigure}{0.49\textwidth}
         \centering
         \includegraphics[width=0.75\textwidth]{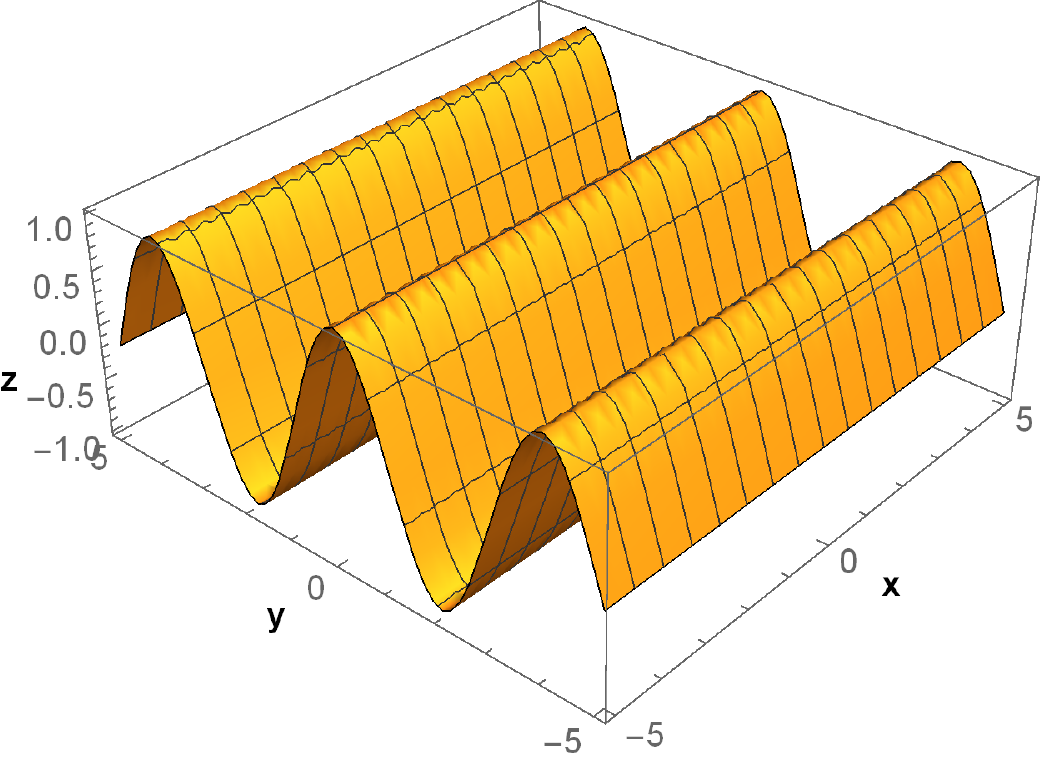}
         \caption{$\cos(ay)$}
         \label{figure:figure0b}
     \end{subfigure}
     \hfill
		\\
		\begin{subfigure}{0.49\textwidth}
         \centering
         \includegraphics[width=0.75\textwidth]{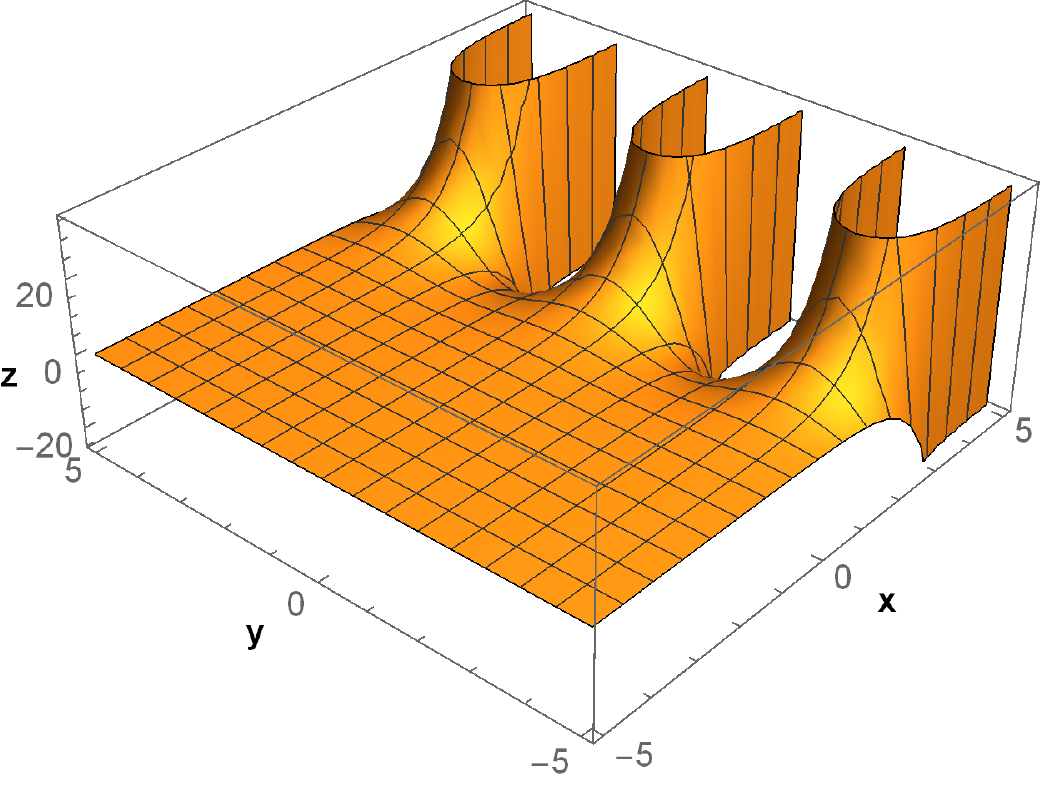}
         \caption{$e^{ax}\cos(ay)$}
         \label{figure:figure0c}
     \end{subfigure}
     \hfill
     \begin{subfigure}{0.49\textwidth}
         \centering
         \includegraphics[width=0.75\textwidth]{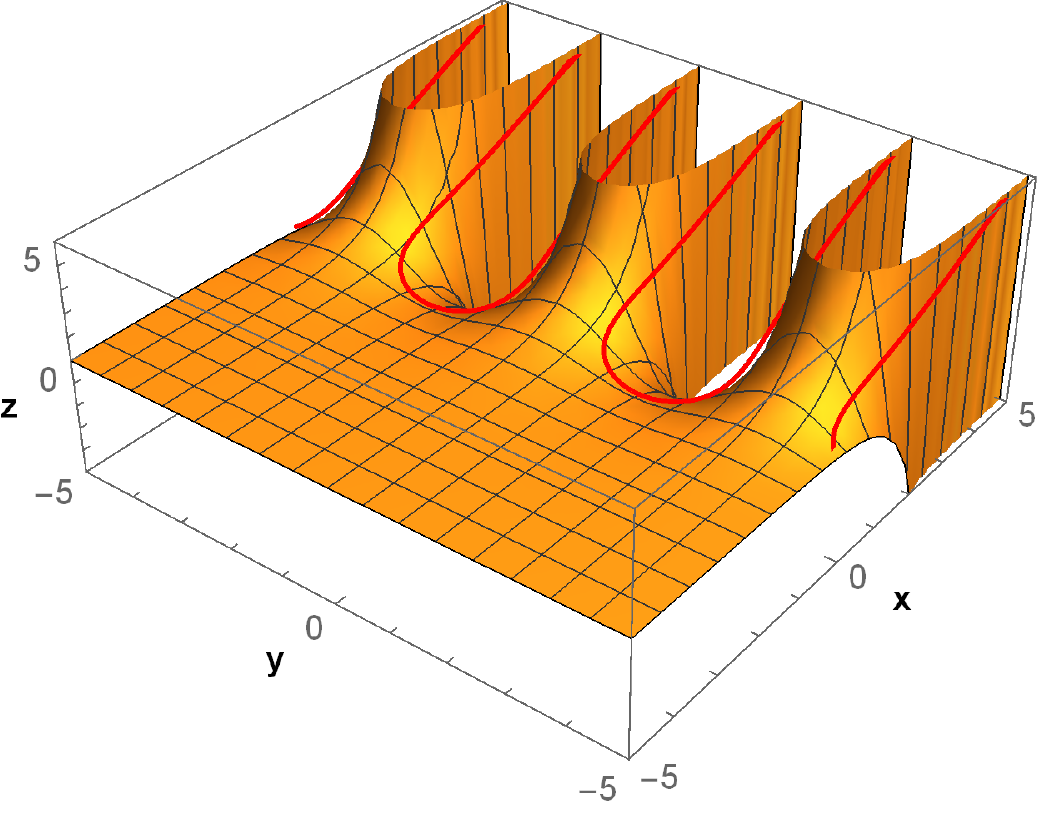}
         \caption{$e^{ax}\cos(ay)=x$}
         \label{figure:figure0d}
     \end{subfigure}
     \hfill
		
     \caption{Creation of contour plots}
        \label{figure:figure0}
\end{figure}

For the $1$-cycle case, in order for $w=z^w$ with $\Log(z)=a+bi$ and $w=x+iy$, the real and imaginary parts of both sides must be equal:
\begin{equation}
\begin{aligned}
x=e^{ax-by}\cos(ay+bx), \\
y=e^{ax-by}\sin(ay+bx).
\end{aligned}
\end{equation}
First consider $z$ a real number.  Then
\begin{equation}
\begin{aligned}
x=e^{ax}\cos(ay),\\
y=e^{ax}\sin(ay).
\label{eqn:eqn1000}
\end{aligned}
\end{equation}
Figure \ref{figure:figure0} shows small component sections of $e^{ax}\cos(ay)$ for $z=2$.  The product of $e^{ax}$ with $\cos(ay)$ creates Figure \ref{figure:figure0c}.  And Figure \ref{figure:figure0d} shows red contour lines over Figure \ref{figure:figure0c} where $e^{ax}\cos(ay)=x$.  A similar diagram can be created for the imaginary expression.  Now consider just the contour lines of (\ref{eqn:eqn1000}) in the $w$-plane.  This is called a contour diagram and used extensively in this paper, an example of which is shown in Figure \ref{figure:figure7}.  In the figure, the real contours are in red and imaginary in blue.  Because $\sin$ and $\cos$ are out of phase, the red and blue contours intersect infinitely often and therefore the zeros of $T_1$ can be enumerated by $\mathbb{Z}$.

\begin{figure}[!ht]
	\centering
			\includegraphics[scale=.8]{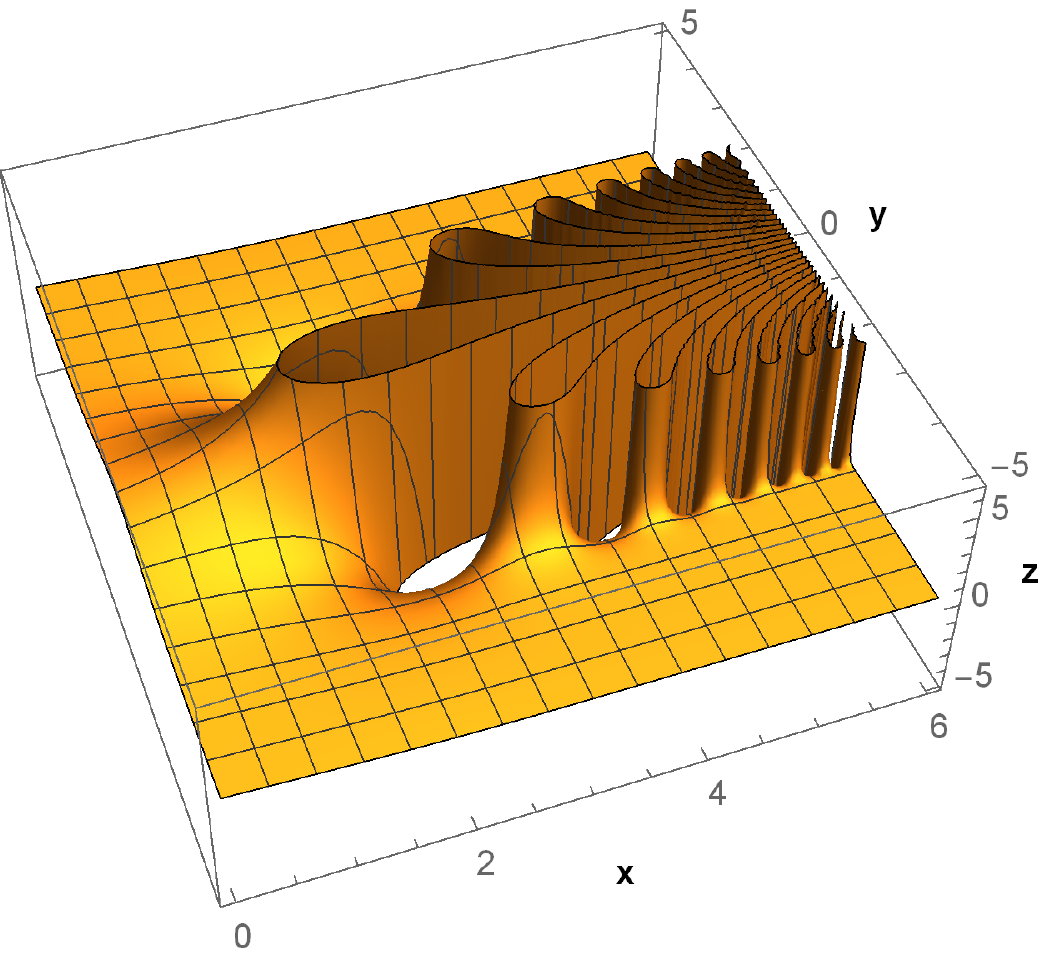}
		\caption{A foliation of a $1$-cycle branch}
   \label{figure:figure1e}
\end{figure}

The analogous $T_2$ expressions
\begin{equation}
\begin{aligned}
y=\exp\{e^{ax}[a\cos(ay)-b\sin(ay)]\}\sin\left(e^{ax}(a\sin(ay)+b\cos(ay))\right)=g(x,y),\\
x=\exp\{e^{ax}[a\cos(ay)-b\sin(ay)]\}\cos\left(e^{ax}(a\sin(ay)+b\cos(ay))\right)=h(x,y)\\
\label{eqn:eqn1001}
\end{aligned}
\end{equation}
are still products of exponential and trigonometric functions although the product is more complicated.  Analogous 3D plots of $T_2$ are foliations of the $T_1$ plots .  Figure \ref{figure:figure1e} shows a small section of one branch or \quotes{contour bulb} of  $h(x,y)$ foliating into \quotes{leaves}. It is reasonable to conjecture this foliation continues indefinitely to the right of the diagram for each branch due to the infinite cycling of the trigonometric functions.   And similar to the $1$-cycle case, contours $x=h(x,y)$ and $y=g(x,y)$ of $T_2$ intersect across the foliated branch leaves at zeros of the function.  Plotting these contours in the $w$-plane produces a $T_2$ contour plot.  An example is shown in Figure \ref{figure:figure20}.  

Therefore, a plausible working hypothesis would be all zeros of $T_2$ can be enumerated by $\mathbb{Z}^2$.  However, in order to establish a connection with the underlying complex geometry, a set of triplets $\{n,m,p\}$ is proposed as an index into the zeros with $(n,m)$ enumerated by $\mathbb{Z}^2$ and $p$ constrained to a small set of positive integers.   The remainder of this paper describes how this indexing is constructed and algorithms to numerically compute the corresponding roots to arbitrary precision.

A major focus of this work establishes a connection between the root ID $\{n,m,p\}$ and a unique branch of $T_2$.  Roots near the origin require special treatment so emphasis is placed on these. The remaining roots are comparatively easier to compute.     
%
%
\section{Nomenclature used in this paper}
\begin{enumerate}
\item
The roots of $T_2$ are organized into \quotes{branch} groups with each root of a branch assigned a \quotes{leaf} number.  This grouping is a reflection of the underlying  logarithmic geometry of $T_2$.  The roots are identified by $\{n,m,p\}$ with $m$ referring to the branch and $n$, the leaf.  The third parameter,$p$, is used when multiple roots are assigned to an $\{n,m\}$ pair and omitted otherwise.
\item
$T_2$ can be placed in "normal" or sometimes called \quotes{rotated} form.   The un-normalized form is called the "base" form.  In rotated form, the branch lobes are horizontal and opening to the right.  Normalization is done to simplify calculations.
\item
Two definitions of "branch" are used:  A contour branch is a group of roots organized on one contour bulb (described above).  A $\pLog$ branch is a single-valued analytic sheet of the composite log function $\pLog$.
\item
The notation $\log(z)$ is the multivalued logarithmic function base $e$.  $\Log(z)$ is the principal value of $\log(z)$.  Likewise, $\arg(z)$ is the multivalued argument function, and $\Arg(z)$ is the principal value of $\arg(z)$.
\item
The composite log function used to compute the roots of $T_2$ has the form of $\log(\log(f(z,w)))$ and is named $\pLog$.   Derived versions of this function are named pLogXXY where \quotes{XX} is a region in the $w$-plane, and \quotes{Y} is the branch-cut type.  For example, $\texttt{pLog3AP}$ refers to Region 3A with branch-cut along the (P)ositive real axis. \texttt{pLog3AN} refers to Region 3A with branch-cut on the negative real axis.  A plot of the real or imaginary component of $\pLog$ is called a  "log stack".  One or more roots are located on single-valued sheets of the log stack having branch number $m$ and leaf number $n$ and correspond to a specific branch and leaf location in a contour plot.  Sequential leaf indexes $n$ for constant $m$ correspond to either a clockwise or counter clockwise rotation around the contour branch leaves.  Thus, there is a close connection between the root ID, a branch of $T_2$, and a contour plot.  Multiple roots on a log sheet are identified by the third index $p$.  In the cases studied in this work, a maximum of $3$ roots were found on a single sheet.  Three stack sheets are often analyzed in this paper and color-coded red for $\{0,1\}$, green for $\{0,0\}$ and blue for $\{-1,0\}$.  
\item
In a case where a root $\{n,m,p\}$ is shown as a point on a plot, root with positive $n$ are color-coded black, and negative $n$, yellow except where noted.
\item
Two auxiliary functions, $\branchF(v,\psi)$ and $\leafF(v,\phi)$ are used.  Plots of these functions are shown in a contour diagram as brown and black traces respectively.
\item
In this paper, precision is the overall number of accurate digits in a result.  Accuracy is the number of digits to the right of the decimal place. Accuracy of a root is determined by back-substitution into the composite log form of $T_2$. This accuracy will be greater than a back-substitution into the exponential form $w-z^{z^w}$.  Using the exponential form to test roots with very large values can cause overflow or underflow.  Therefore, the logarithmic form is used to check all roots.   
\item
The roots are computed by fixed-point iteration of the composite log expression $w-\pLog(w,n,m;z)$ or one of it's derived forms.  For each value of $m$, a seed is computed close to the head of the contour bulb.  Seed points in a graph are green.  The seed is initially computed with arbitrary precision and set to a desired numeric precision for the iteration.  In order to achieve a desired accuracy, the seed precision is set to a value higher than the desired accuracy.   In the case of computing a sequential set of roots on the same branch bulb, the seed for the next root is set to the previous root at a precision of the initial seed.   This is because sequential roots on the same branch bulb grow closer together as $n$ is incremented (or decremented in the negative case).  
\end{enumerate}

A straight-forward Newton iteration of $w-z^{z^w}$ is problematic:  The quantity $z^{z^w}$ is interpreted by most software packages as principal-valued and the exponential when iterated often leads to underflow or overflow. However, if the function is converted to its logarithmic form, the root basins for a majority of roots tested encompass a large area in the vicinity of the roots, and the logarithmic terms are must less susceptible to underflow and overflow.  Given $w=z^{z^w}$, then
\begin{equation}
\begin{aligned}
w&=e^{\Log(z)e^{w\Log(z)}}\\
(\Log(w)+2n\pi i)&=\Log(z)e^{w\Log(z)};\quad n\in\mathbb{Z} \\
\frac{1}{\Log(z)}\bigg[\Log(w)+2n\pi i\bigg]&=e^{w\Log(z)}\\
w&=\frac{1}{\Log(z)}\bigg\{\Log\bigg[\frac{1}{\Log(z)}\left(\Log(w)+2n\pi i\right)\bigg]+2m\pi i\bigg\};\quad (n,m)\in\mathbb{Z}\\
w&=\pLog(w,n,m;z)
\end{aligned}
\label{equation:eqn1}
\end{equation}
so that $T_2$ in it's logarithmic form is
\begin{equation}
T_2(w,n,m;z)=w-\pLog(w,n,m;z);\quad (n,m)\in\mathbb{Z}.
\end{equation}
This particular form of $T_2$ offers the advantage of easily generating roots with very high values of $n$ and $m$ although some precautions are needed to adjust the method for branch-cuts.  The function $\pLog$ is a nested logarithm with a nested logarithmic geometry.  To better understand this geometry, first consider $\log(z)=\Log(|z|)+i\arg(z)$.  In the complex plane, a plot of $\arg(z)$, is a helical coil winding around the origin an infinite number of times.  The coil can be partitioned into single-valued analytic sheets or \quotes{branches} indexed by an integer $k$ according to the definition of the complex logarithm: $\log(z)=\Log(z)+2k\pi i$.

Both the real and imaginary components of $\pLog$ have a primary "m" coil, with each single-valued sheet of the m-coil analytically-continuing in both clockwise and counter clockwise directions to secondary "n" coils winding around the origin an infinite number of times in both directions.  This branching can be partitioned into single-valued sheets or \quotes{leaves} indexed by the pair $\{n,m\}$.  One or more roots $\{n,m,p\}$ of $T_2$ will be located on sheet $\{n,m\}$. 
\begin{figure}[!ht]
	\centering
			\includegraphics[scale=.5]{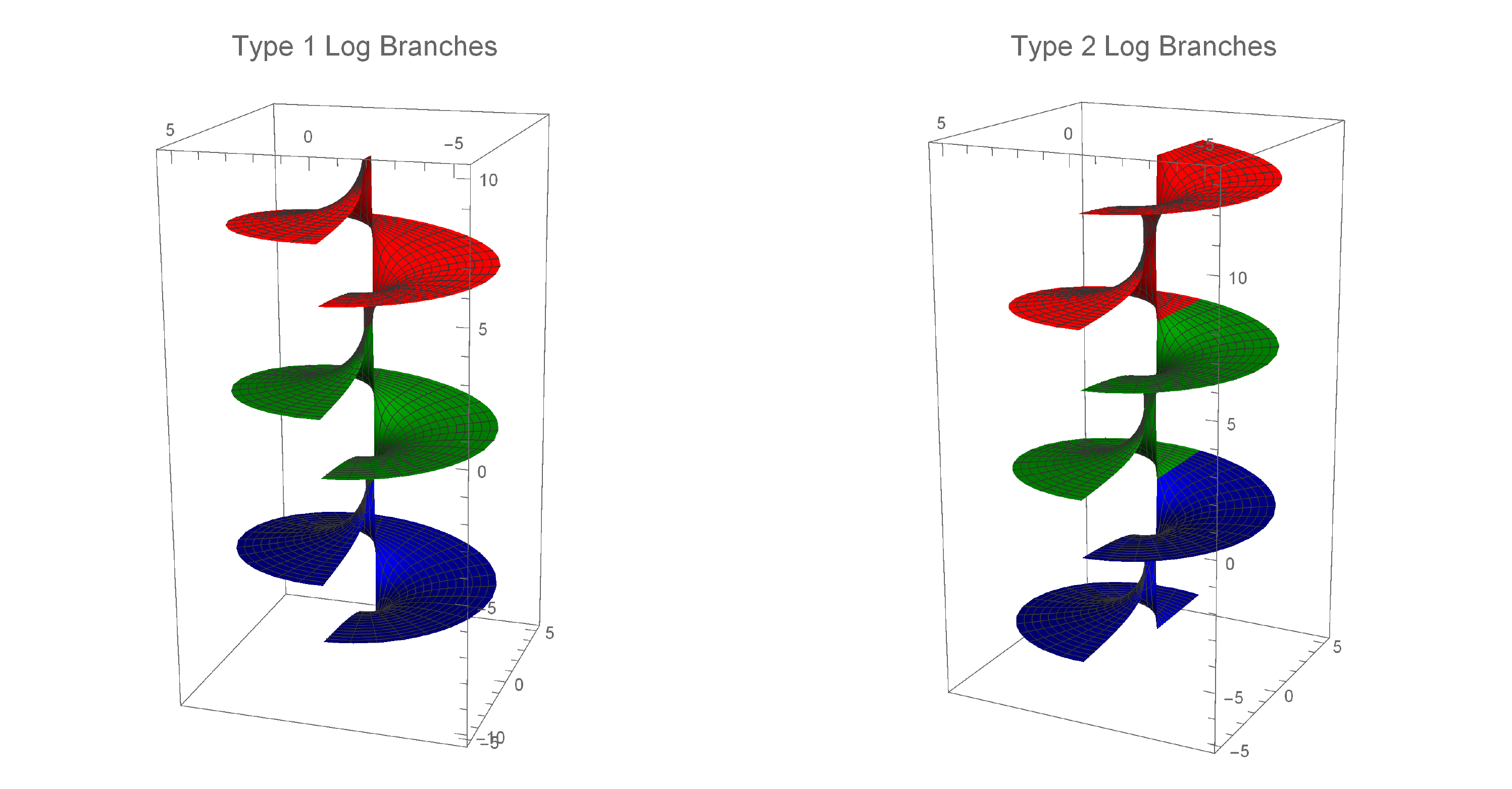}
		\caption{Three primary sheets of $\log(\log(z))$}
   \label{figure:figure1}
\end{figure}      

First consider the simple case of log(log(z)).  It has two singular points and as a consequence, single-valued analytic sheets can be defined in two ways: 
\begin{enumerate}
\item 
Type-I Log Branches. Branch-cut: $(-\infty,1)$.
\item
Type-II Log Branches.  Branch-cuts: $ (-\infty,0) \cup (1,\infty)$.
\end{enumerate}

The root finding algorithms described below splice together sheets of $\pLog$ as necessary to provide a contiguous and analytic surface to iterate over.   For example,  consider the first plot in Figure \ref{figure:figure1} showing three primary sheets of $\log(\log(z))$.  The green sheet is $\{0,0\}$, the red is $\{0,1\}$, and blue is $\{0,-1\}$.  This is a Type-I log stack with a single branch cut $ (-\infty,1)$.  However, the half-section of the green sheet in the upper half $w$-plane can be spliced with the half-section of the red sheet in the lower half plane.  Doing likewise with the other sheets create a Type-II log stack shown in the second diagram. 

It's important to remember Figure \ref{figure:figure1} depicts only three of an infinite number of branches in the primary coil.  And each single-valued primary branch in turn analytically-continues in both directions to secondary coils each with an infinite number of leaf branches. Understanding this geometry and its connection to the root indexing $\{n,m,p\}$ is central to understanding the algorithms described in this paper. 

The complex w-plane is partitioned into four main regions:
\begin{enumerate}
\item
 Region 1: the positive real line,
\item
 Region 2: inside the unit circle,
\item
 Region 3:  the unit circle,
\item
 Region 4: outside the circle.
\end{enumerate}

Region 1 is further sub-divided into six sub-regions 1A through 1F.  Regions 2 through 4 are sub-divided into parts A and B.   Each region has an associated $\pLog$ iterator.  

In order to more clearly present the transformations used for the 2-cycle case, they are first applied to the 1-cycle case.
%
%
\section{1-cycle iterated exponential}

For $T_1(w;z)=w-z^w$, let $\Log(z)=a+bi$ and $w=x+iy$.  Then for $T_1=0$,
\begin{equation}
\begin{aligned}
y&=e^{ax}\sin(ay),\\
x&=e^{ax}\cos(ay).
\end{aligned}
\label{eqn:eqn1020}
\end{equation}

\begin{figure}[!ht]
     \centering
     \begin{subfigure}{0.4\textwidth}
         \centering
         \includegraphics[width=0.75\textwidth]{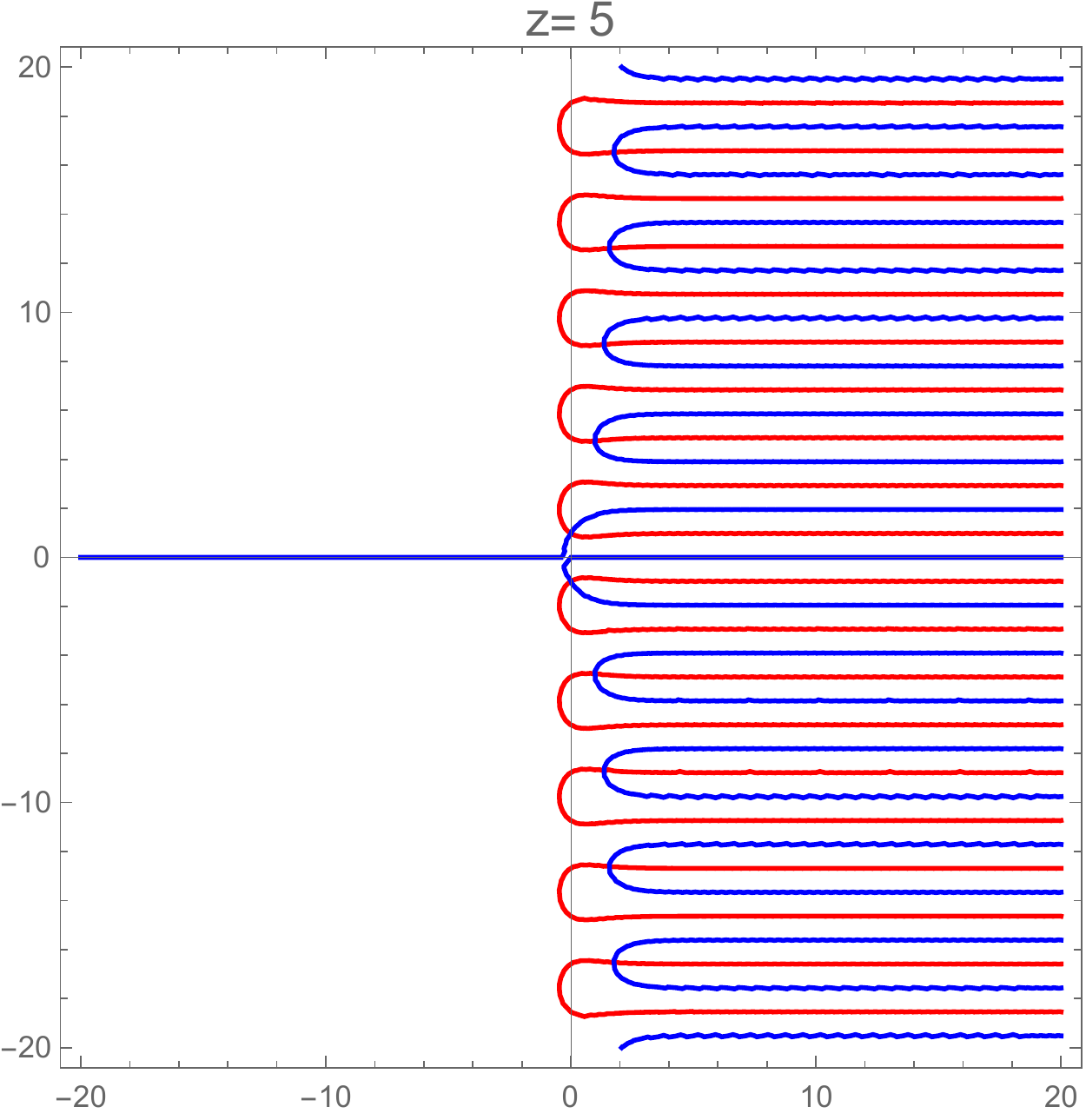}
         \caption{1-cycle contours}
         \label{figure:figure7a}
     \end{subfigure}
     \hfill
     \begin{subfigure}{0.4\textwidth}
         \centering
         \includegraphics[width=0.75\textwidth]{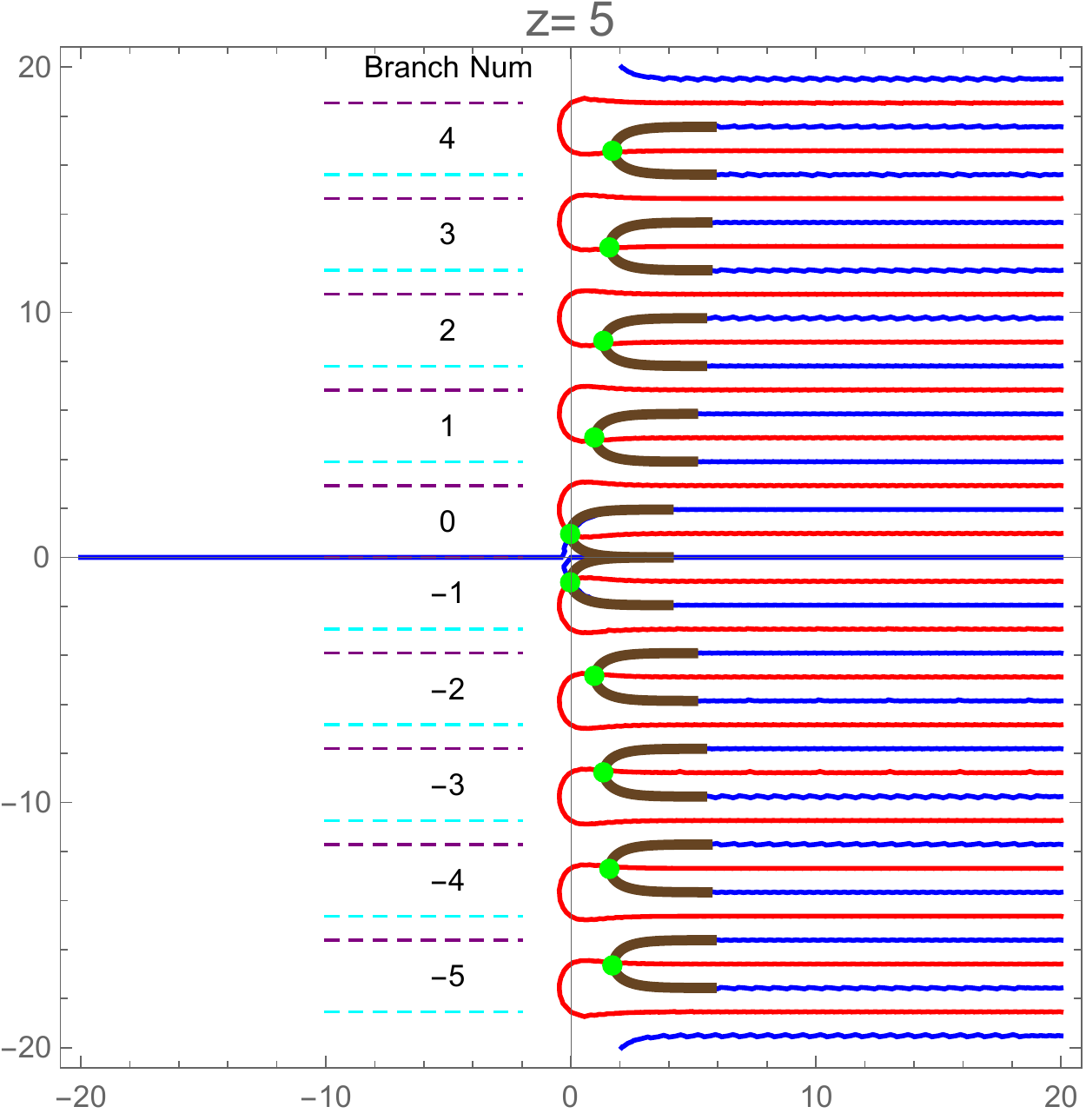}
         \caption{1-cycle contours with branchF[v,$\psi$] traces}
         \label{figure:figure7b}
     \end{subfigure}
     \hfill
     \caption{1-cycle contour plots}
        \label{figure:figure7}
\end{figure}

$T_1$ contour plots are shown in Figure \ref{figure:figure7}.    The branching is already in normal form.  

Set the right sides of (\ref{eqn:eqn1020})  to real constants:

\begin{equation}
\begin{aligned}
e^{ax} \sin(ay)&=\psi,\\
e^{ax} \cos(ay)&=\phi.
\end{aligned}
\label{equation:eqn4c}
\end{equation}

Solving for $x$ in both cases define:

\begin{equation}
\begin{aligned}
\branchF(y,\psi)&=1/a \Log[\frac{\psi}{\sin{ay}}],\\
\leafF(y,\phi)&=1/a \Log[\frac{\phi}{\cos{ay}}].
\end{aligned}
\label{equation:eqn4d}
\end{equation}

The asymptotes of (\ref{equation:eqn4d}) are 
\begin{equation}
\begin{aligned}
\text{branchA}(k_1)&=\frac{k_1\pi}{a},\\
\text{leafA}(k_2)&=\frac{\pi}{2a}+\frac{k_2\pi}{a};\quad k\in\mathbb{Z}.
\end{aligned}
\end{equation}
Since the distance between asymptotes of both functions is $\pi/a$, the mean between each asymptote of $\branchF$ is
\begin{equation}
\text{branchM}(k)=\frac{\pi}{2a}+\frac{k\pi}{a}.
\end{equation}

The width of each of the red and blue contour lobes in Figure \ref{figure:figure7} is $\pi/a$.  Thus a \quotes{branch size} can be defined as the combined width across a set of intersecting red and blue lobes.  This is equal to $\frac{3\pi}{2a}$.  The branch number becomes $k$.  
Now consider the blue components of the contour plot in the upper half-plane.  These are equal to the value of the imaginary part, $y$ of $w$.  In order to obtain an approximation to these contours, set $\psi$ equal to the median of the blue contours which is $\text{branchM}(k)$.  

Figure \ref{figure:figure7b} displays small segments of $\branchF$ as the brown contours and note how closely they follow the blue contours. If the analytic expression for the roots of this expression was unknown, as in the case of $T_2$, this would provide an easy method to compute iteration seeds close to the roots:
\begin{equation*}
\text{seedF}(k)=\{\branchF(\psi_k,\psi_k),\psi_k\},\quad \psi_k=\text{branchM}(k).
\end{equation*}
\section{1-cycle case when z is not a positive number greater than one}

\begin{figure}[!ht]
     \centering
     \begin{subfigure}{0.4\textwidth}
         \centering
         \includegraphics[width=0.75\textwidth]{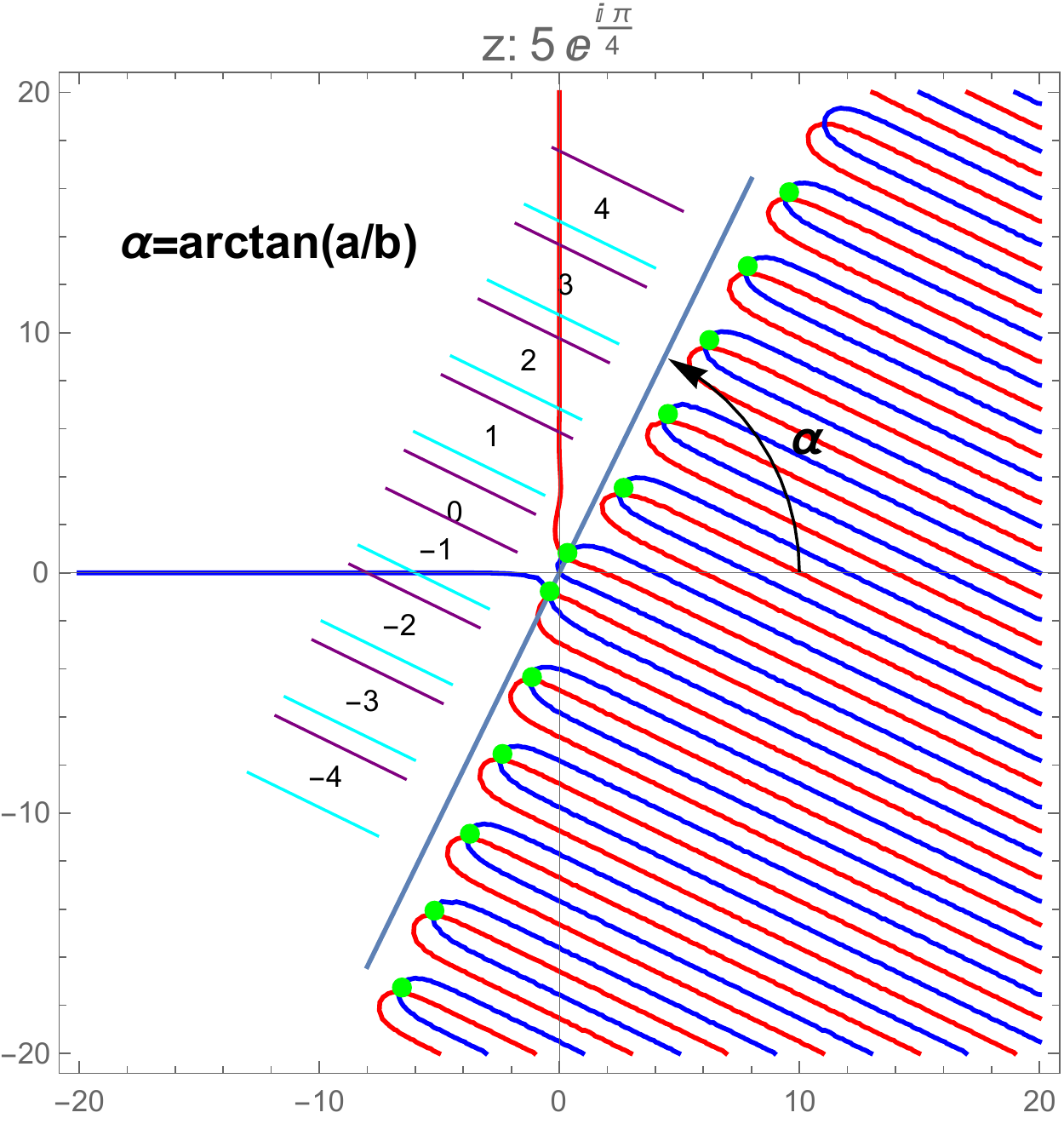}
         \caption{1-cycle contours in base form}
         \label{figure:figure6a}
     \end{subfigure}
     \hfill
     \begin{subfigure}{0.4\textwidth}
         \centering
         \includegraphics[width=0.75\textwidth]{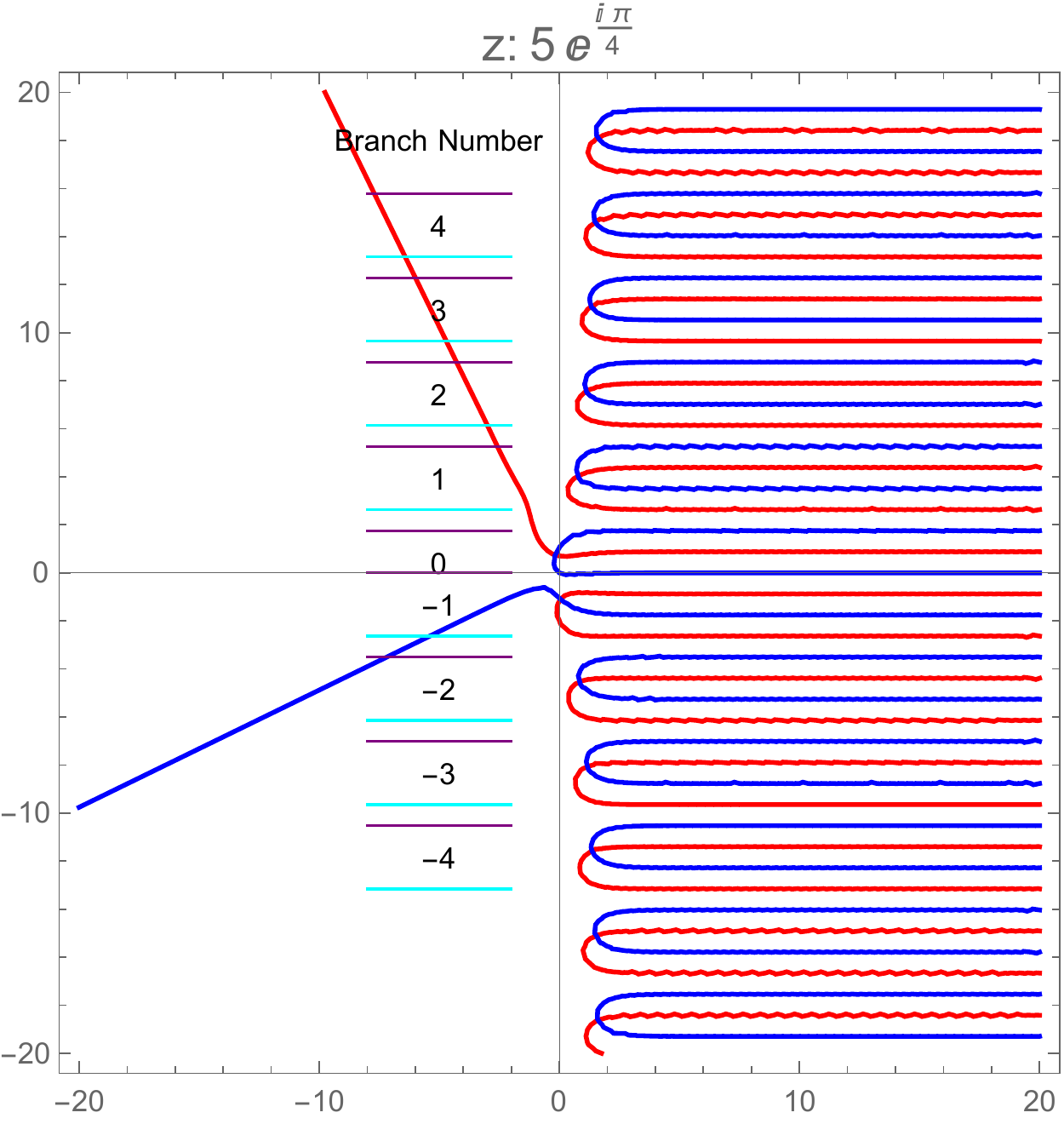}
         \caption{1-cycle contours in normal form}
         \label{figure:figure6b}
     \end{subfigure}
     \hfill
     \caption{1-cycle contour plot types}
        \label{figure:figure6}
\end{figure}

In the cases when $z$ is not a positive real greater than 1, the real and imaginary contour lobes are no longer horizontal but are tilted with respect to the axes.  For example, $z=5 e^{\pi i/4}$ produces the contour plots in Figure \ref{figure:figure6}. 

Notice  the contour branching in base form is tilted with respect to the coordinate axes.  The angle of inclination, $\alpha$, in Figure \ref{figure:figure6a} is easily calculated.  Letting $\Log(z)=\ln(r)+i\theta=a+bi$ and $w=x+iy$, then
\begin{equation}
\begin{aligned}
z^w&=e^{(a+bi)(x+iy)}=e^{(ax-by)+i(ay+bx)}\\
&=e^{(ax-by)}\big[\cos(ay+bx)+i\sin(ay+bx)\big].
\end{aligned}
\label{equation:eqn12}
\end{equation}
And for $w=z^w$, define the (B)ase forms of the contours as
\begin{equation}
\begin{aligned}
\text{imagF1B}(x,y)&=e^{(ax-by)}\sin(ay+bx)=y,\\
\text{realF1B}(x,y)&=e^{(ax-by)}\cos(ay+bx)=x.
\label{equation:eqn15}
\end{aligned}
\end{equation}

The sine term is maximum when $\displaystyle ay+bx=\frac{\pi}{2}+2k\pi$ or $\displaystyle y=-\frac{b}{a} x+\frac{\pi(2k+1/2)}{a}$.  The negative reciprocal of the slope of this line is the slope of the contours.  Thus the contours make an angle of $\alpha=\arctan(a/b)$.  

Placing the contours in normal form with the lobes horizontal and opening to the right simplifies the the analysis.  The simplest way to do this is to rotate the contours in Figure \ref{figure:figure6a} by an angle of $\pi/2-\alpha$.  Letting $\beta=\pi/2-\alpha$, define

\begin{equation}
\text{rotationF}(\{x,y\},-\beta)=\big\{x\cos(\beta)+y\sin(\beta),y\cos(\beta)-x\sin(\beta)\big\}=\{u,v\}.
\label{equation:eqn15b}
\end{equation}
Composing (\ref{equation:eqn15}) with (\ref{equation:eqn15b}), we obtain the (N)ormal forms
\begin{equation}
\begin{aligned}
\text{imagF1N}(x,y)&=\text{imagF1B}\big(u,v)\big)\equiv v,\\
\text{realF1N}(x,y)&=\text{realF1B}\big(u,v)\big)\equiv u.
\label{equation:eqn16}
\end{aligned}
\end{equation}

A contour plot of (\ref{equation:eqn16}) is shown in Figure \ref{figure:figure6b}. 

Once the equations are in normal form, branch sizes, auxiliary functions, optimum seed locations and other parameters are easily calculated.  For example, in order to compute the $T_1$ auxiliary functions, the variables would need to be separated in the expressions 
\begin{equation}
\begin{aligned}
e^{(ax-by)}\sin(ay+bx)&=\psi, \\
e^{(ax-by)}\cos(ay+bx)&=\phi.
\label{equation:eqn16b}
\end{aligned}
\end{equation}
But this cannot be done in closed.  However, the variables in normal form can be separated.  That is, if we have
\begin{equation}
\begin{aligned}
\text{imagF1B}\big(u,v\big)=\psi, \\
\text{realF1B}\big(u,v\big)=\phi,
\label{equation:eqn17}
\end{aligned}
\end{equation}
the variables $(x,y)$ are separable and will turn out to be so in the 2-cycle case.  Solving both expressions in (\ref{equation:eqn17}) for $x$ in terms of $y$, define 
\begin{equation}
\begin{aligned}
\branchF(y,\psi)=\frac{4\Log\bigg[\psi\csc(1/4 y\sqrt{\pi^2+16\Log^2(5)})\bigg]}{\sqrt{\pi^2+16\Log^2(5)}},\\
\leafF(y,\phi)=\frac{4\Log\bigg[\phi\sec(1/4 y\sqrt{\pi^2+16\Log^2(5)})\bigg]}{\sqrt{\pi^2+16\Log^2(5)}}.\\
\end{aligned}
\end{equation}
The asymptotes are solutions to
\begin{equation}
\begin{aligned}
\sin\bigg[1/4 y\sqrt{\pi^2+16\Log^2(5)}\bigg]&=0, \\
\cos\bigg[1/4 y\sqrt{\pi^2+16\Log^2(5)}\bigg]&=0.
\label{eqn:eqn1b2}
\end{aligned}
\end{equation}

Solving (\ref{eqn:eqn1b2}) for the rotated dimensions of the real and imaginary asymptotes
:
\begin{equation}
\begin{aligned}
\text{imag contours }&: \biggr\{\frac{8k\pi}{\sqrt{\pi^2+16\Log^2(5)}},\frac{4(\pi+2k\pi)}{\sqrt{\pi^2+16\Log^2(5)}}\biggr\},\\
\text{real contours }&: \biggr\{\frac{4(-\pi/2+2k \pi)}{\sqrt{\pi^2+16\Log^2(5)}},\frac{4(\pi/2+2k \pi)}{\sqrt{\pi^2+16\Log^2(5)}}\biggr\}.
\label{equation:eqn17b}
\end{aligned}
\end{equation}
From (\ref{equation:eqn17b}), the size (width) of each blue and red contour bulb is
$$
\epsilon=\frac{4\pi}{\sqrt{\pi^2+16\Log^2(5)}},
$$
so that the total branch size in this case (from purple to cyan lines) is:
$$
\delta=\frac{6\pi}{\sqrt{\pi^2+16\Log^2(5)}}.
$$
However, all values computed are relative to the rotated frame.  In order to compute branch parameters such as seeds, the data is inverted back into the base frame by inverting the rotation transformation.  For example, given the rotated seed $(x,y)$, invert it via $\text{rotationF}[\{x,y\},\beta]$,  to obtain the real and imaginary components of the base seed:
\begin{equation}
\text{baseSeed=}\big\{x\cos(-\beta)+y\sin(-\beta),(y\cos(-\beta)-x\sin(-\beta))\big\}.
\end{equation}

Figure \ref{figure:figure6a} shows the completed base-frame contour plot for this example.  A similar procedure is followed for the 2-cycle case.
%
%
\section{2-cycle cases}

The techniques used to analyze the 1-cycle case are now applied to the 2-cycle case by first  deriving 2-cycle expressions for the auxiliary functions. 

\begin{figure}
     \centering
     \begin{subfigure}{0.4\textwidth}
         \centering
         \includegraphics[width=0.75\textwidth]{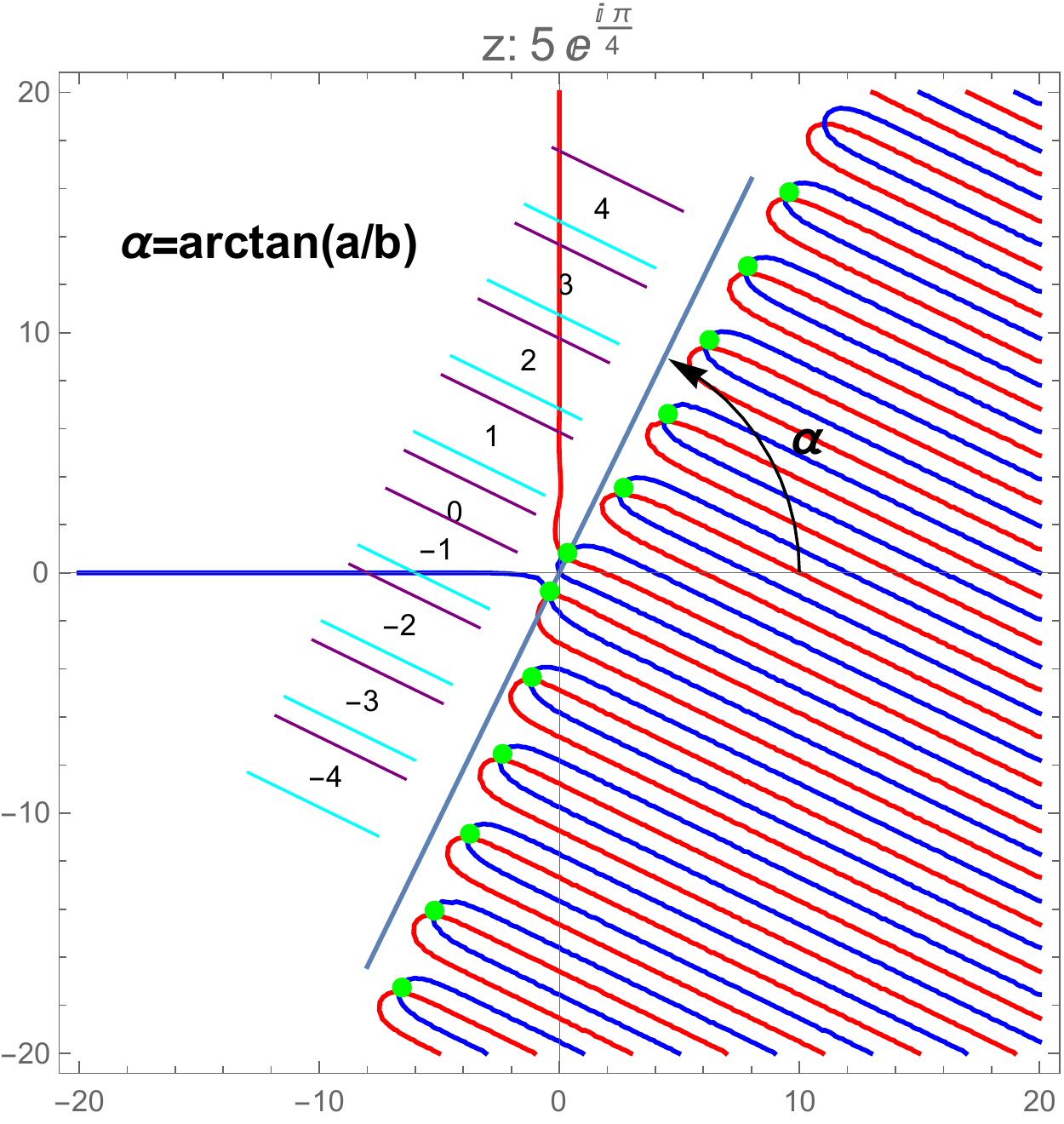}
         \caption{1-Cycle contour plot}
         \label{figure:figure12a}
     \end{subfigure}
     \hfill
     \begin{subfigure}{0.4\textwidth}
         \centering
         \includegraphics[width=0.75\textwidth]{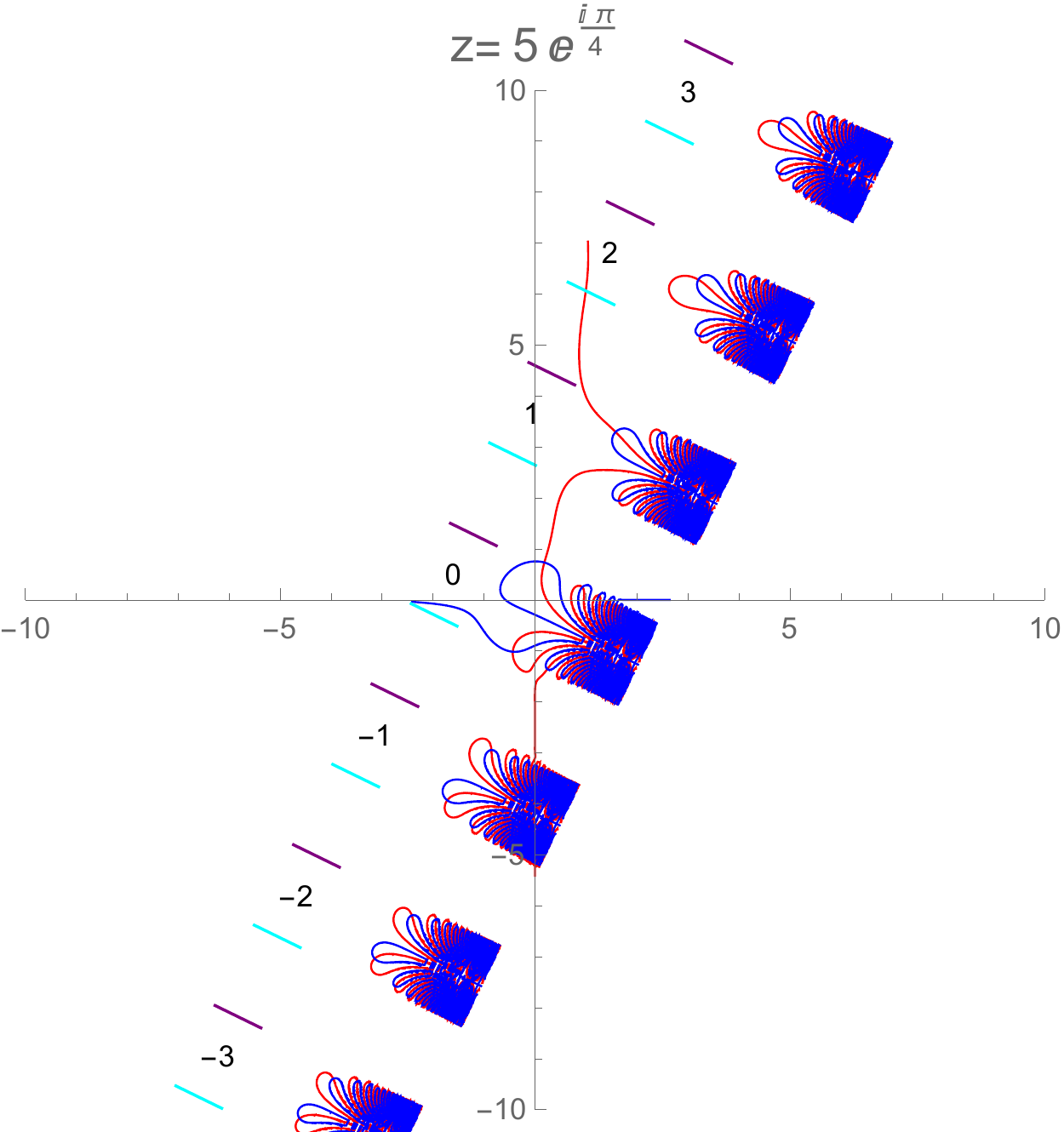}
         \caption{2-cycle contour plot}
         \label{figure:figure12b}
     \end{subfigure}
     \hfill
     \caption{Comparison of cycle types}
        \label{figure:figure12}
\end{figure}

Consider the contour plots for a 1-cycle:   A 2-cycle iterated exponential  can be considered a composition of those branches by virtue of the expression $z^{z^w}$.  This leads to a foliation of each 1-cycle branch into an infinite set of 2-cycle leaves.  Each 2-cycle leaf, including degenerate ones near the origin, can be identified by a root ID $\{n,m,p\}$.  An example of this foliation is shown in Figure \ref{figure:figure12} which compares  the $T_1$ contour with the $T_2$ contour diagram for $z=5 e^{\pi i/4}$.  As with the 1-cycle, the real plots are in red and the imaginary are in blue.  The 2-cycle branches are delimited by purple and cyan lines with the branch numbers derived from the asymptote expressions of $T_2$.  Each intersection of a red and blue leaf represents a root of $T_2$ although some intersections near the origin are not in the form of bulbs.   $T_2$ is first normalize  to compute branch parameters, auxiliary equations, and normalized seeds.  The normalized seeds are then inverted as starting values for a fixed-point iteration of branch roots. 

In the case of $T_2$, the expressions for the real and imaginary contours are more complicated.  Recall the $1$-cycle expressions:
\begin{equation}
\begin{aligned}
z^w&=e^{(a+bi)(x+iy)}=e^{(ax-by)+i(ay+bx)}\\
&=e^{(ax-by)}\bigg[\cos(ay+bx)+i\sin(ay+bx)\bigg].
\end{aligned}
\label{equation:eqn21}
\end{equation}
Now let
\begin{equation}
\begin{aligned}
e^{(ax-by)}\cos(ay+bx)&=c \\
e^{(ax-by)}\sin(ay+bx)&=d.
\label{equation:eqn22}
\end{aligned}
\end{equation}
Then
$$
\begin{aligned}
z^{z^w}&=e^{s(a+bi)}=e^{(a+bi)(c+di)}\\ 
&=\text{exp}\big\{(ac-bd)\left[\cos(ad+bc)+i\sin(ad+bc)\right]\big\}=x+iy,
\end{aligned}
$$
so that in order for $w=z^{z^w}$:
$$
\begin{aligned}
y&=e^{(ac-bd)}\sin(ad+bc)\\
x&=e^{(ac-bd)}\cos(ad+bc).\\
\label{equation:eqn22b}
\end{aligned}
$$

Plugging in $c$ and $d$ above, define

\begin{equation}
\begin{split}
\text{imagF2B}(x,y)&=\text{exp}\bigg\{e^{ax-by}\big[a\cos(ay+bx)-b\sin(ay+bx)\big]\bigg\}\\
&\hspace{20pt}\times\sin\bigg\{e^{ax-by}\big[a\sin(ay+bx)+b\cos(ay+bx)\big]\bigg\},
\end{split}
\end{equation}

\begin{equation}
\begin{split}
\text{realF2B}(x,y)&=\text{exp}\bigg\{e^{ax-by}\big[a\cos(ay+bx)-b\sin(ay+bx)\big]\bigg\}\\
&\hspace{20pt}\times\cos\bigg\{e^{ax-by}\big[a\sin(ay+bx)+b\cos(ay+bx)\big]\bigg\}.
\end{split}
\end{equation}

Now assign the arguments of the outer exp and trig terms to constants:

\begin{equation}
\begin{aligned}
e^{ax-by}\big[a\cos(ay+bx)-b\sin(ay+bx)\big]&=\psi, \\
e^{ax-by}\big[a\sin(ay+bx)+b\cos(ay+bx)\big]&=\phi.
\end{aligned}
\label{eqn:eqn100}
\end{equation}

As written, each equation of (\ref{eqn:eqn100}) cannot be solved for $x$ in terms of $y$ in closed form.  However the variables are separated in the rotated forms.  For example, letting $z=5 e^{\pi i/4}$ and applying the rotation transformation to the left sides of (\ref{eqn:eqn100}), gives:  
\begin{equation}
\begin{split}
&-1/4e^{1/4x\sqrt{\pi^2+16\Log^2(5)}}\times \\
&\hspace{20pt}\bigg(-4\cos\big[1/4 y\sqrt{\pi^2+16\Log^2(5)}\big]\Log(5)+\pi\sin\big[1/4 y\sqrt{\pi^2+16\Log^2(5)}\big]\bigg)=\psi,\\
\\
&1/4e^{1/4x\sqrt{\pi^2+16\Log^2(5)}}\times \\
&\hspace{20pt}\bigg(\pi\cos\big[1/4 y\sqrt{\pi^2+16\Log^2(5)}\big]+4\Log(5)\sin\big[1/4 y\sqrt{\pi^2+16\Log^2(5)}\big]\bigg)=\phi,
\end{split}
\end{equation}
and note how the variables have become separated.  Solving for $x$ in terms of $y$, define:
\begin{equation*}
\begin{aligned}
\branchF(y,\psi)&=\frac{4\Log\big[-\frac{4\psi}{-4\cos\big[1/4 v\sqrt{\pi^2+16\Log^2(5)}\big]\Log(5)+\pi\sin\big[1/4 v\sqrt{\pi^2+16\Log^2(5)}\big]}\big]}{\sqrt{\pi^2+16\Log^2(5)}},\\
\leafF(y,\phi)&=\frac{4 \Log \left[\frac{4 \phi }{\pi  \cos \left(\frac{1}{4} y \sqrt{\pi ^2+16 \log ^2(5)}\right)+4 \log (5) \sin \left(\frac{1}{4} y \sqrt{\pi ^2+16 \log ^2(5)}\right)}\right]}{\sqrt{\pi ^2+16 \log ^2(5)}}.
\end{aligned}
\end{equation*}

The asymptotes of $\branchF$ will determine the branch numbers $m$.   And $\leafF(v,n)$ will trace approximately, the leaves of each branch.  Next, define the asymptote expressions:

\begin{equation}
\begin{aligned}
\text{branchA}(m)&=\left\{\frac{4 \left(2 \pi  m+\tan ^{-1}\left(\frac{4 \log (5)}{\pi }\right)\right)}{\sqrt{\pi ^2+16 \log ^2(5)}},\frac{4 \left(2 \pi  k-\pi +\tan ^{-1}\left(\frac{4 \log (5)}{\pi }\right)\right)}{\sqrt{\pi ^2+16 \log ^2(5)}}\right\},\\
\text{leafA}(n)&=\left\{\frac{4 \left(2 \pi  n-\tan ^{-1}\left(\frac{\pi }{4 \log (5)}\right)\right)}{\sqrt{\pi ^2+16 \log ^2(5)}},\frac{4 \left(2 \pi  n+\pi -\tan ^{-1}\left(\frac{\pi }{4 \log (5)}\right)\right)}{\sqrt{\pi ^2+16 \log ^2(5)}}\right\}.
\end{aligned}
\end{equation}

Recall $\branchF$ was derived from the imaginary component of (\ref{equation:eqn22b}) which was equated to the imaginary part of $w=x+iy$.  In normal form, set $\psi=\pm\Log(\pm|m|)$ as necessary to generate a real value of $\branchF$ where $m$ is the mean of the branch dimension, or a small value if $m=0$.  Then $\branchF(v,\psi)$ is a good approximation to an envelope of a branch within its domain of definition.  Likewise, setting $\phi=2n\pi$ is a good approximation to the trace of leaf $n$.  

With this background, the roots of $T_2$ can now be computed.
%
%
\section{Region 1 Analysis}
\begin{figure}[!ht]
	\centering
			\includegraphics[scale=1]{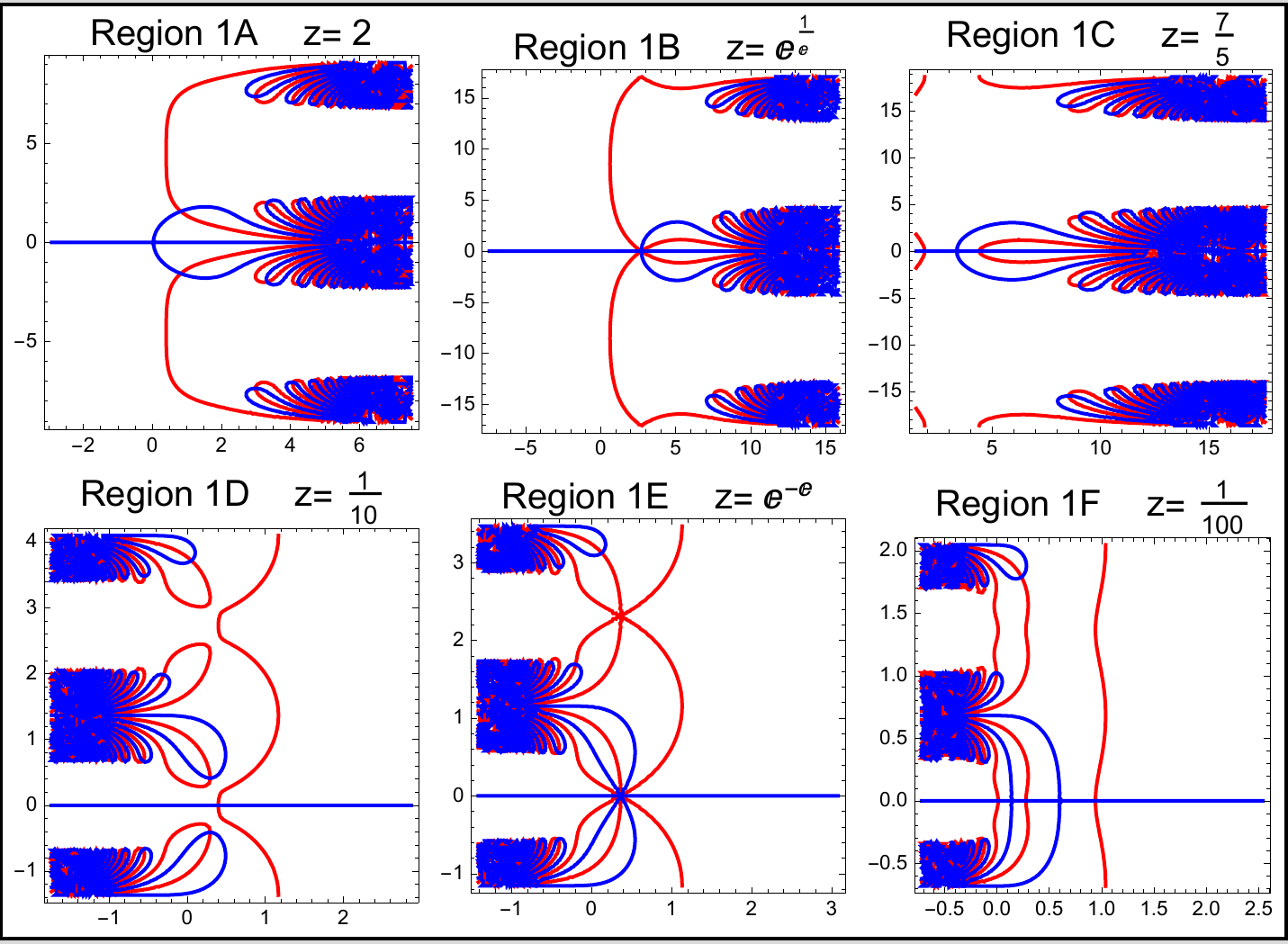}
		\caption{Contour morphologies of Region 1 in base form}
   \label{figure:figure15}
\end{figure}

Region 1 is the positive real axis.  This region is sub-divided into six sub-regions determined by the contour morphologies of branch lobes near the origin.  Example branch contours of each sub-region are shown in Figure \ref{figure:figure15} and described in Table \ref{table:table1}.  Note regions 1A,1B and 1C are already in normal form.

\begin{table}[!ht]
\centering
\begin{tabular}{|c|c|p{12cm}|}
\multicolumn{3}{ c }{} \\
\hline
Region & Domain & Description \\
\hline
 1A & $|z|>e^{1/e}$ &Two leaf zero roots, $\{0,0,1\}$ and $\{0,0,2\}$  \\ 
\hline
 1B & $z=e^{1/e}$ & The two real roots of Region 1A coalesce \newline into a single root at $z=e$ with multiplicity $2$ \\ 
\hline
1C &  $1<|z|<e^{1/e}$&   The single root of Region 1B splits into two real roots\newline $\{0,0,1\}$ and $\{0,0,2\}$\\ 
\hline
1D & $e^{-e}<|z|<1$. & Orientation of the branch leaves invert as the two real\newline roots of Region 1C split into 2 complex roots and one real root,\newline $\{0,0,1\},\{0,0,2\}$ and $\{0,0,3\}$\\ 
\hline
1E &  $|z|=e^{-e}$&   As z approaches $e^{-e}$ from the right, the three roots of \newline Region 1D coalesce into a single root\newline at $e^{-e}$ with multiplicity 3\\ 
\hline
1F &  $0<|z|<e^{-e}$ & The single root of Region 1E splits into three real roots,\newline $\{0,0,1\},\{0,0,2\}$ and $\{0,0,3\}$\\
 \hline
\end{tabular}
\caption{Region $1$ sub-regions}
\label{table:table1}
\end{table}

%
%
\subsection{Region 1A}
This is the simplest case.  A sequence of steps is formulated in this section and used for succeeding sections.

\subsubsection{Analyzing the branch-cuts}
The most important step in computing the roots for a particular region is analyzing the branch-cuts of $\pLog$. Letting $z=2$ as an example in this region, recall the expressions

\begin{equation}
T_2(w,n,m;z)=w-\pLog(w,n,m;z)
\end{equation}
with 
\begin{equation}
\pLog(w,n,m)=\frac{1}{\Log(z)}\bigg\{\Log\bigg[\frac{1}{\Log(z)}\left(\Log(w)+2n\pi i\right)\bigg]+2m\pi i\bigg\}.\\
\end{equation}
The primary branch-cut comes from the $\Log(w)$ term. In order to find the secondary branch-cut, let $w=re^{i\theta}$ and $\Log(z)=a+bi$, then 
\begin{align*}
\Log\big[\frac{1}{a+bi} (\Log[w]+2n\pi i)\big]&= \Log[1/(a+bi) (ln(r)+i\theta+2 n \pi i)] \\
&=Log\bigg[\frac{a\ln(r)+bk}{a^2+b^2}+i\frac{ak-b\ln(r)}{a^2+b^2}\bigg];\quad k=\theta+2 n\pi.
\end{align*}

\begin{figure}
	\centering
			\includegraphics[scale=.5]{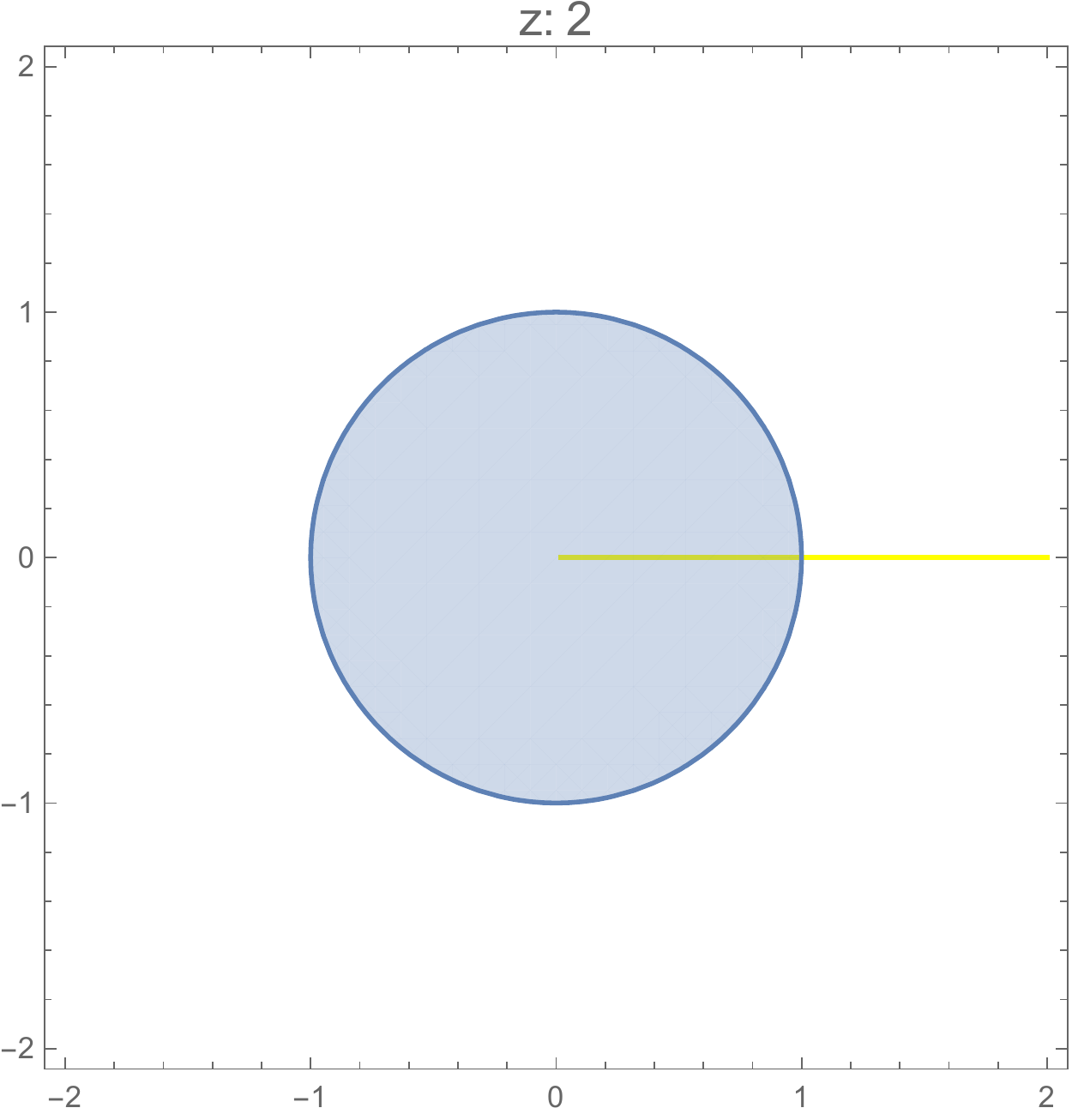}
		\caption{Secondary branch-cut of Region 1A}
   \label{figure:figure16}
\end{figure}

So that the secondary cut occurs when $a\ln(r)+bk<0$ and $ ak-b\ln(r)=0$.  These are the cut-Domain and Trace in this paper. The secondary cut is thus the intersection of these:

\begin{equation}
(a \ln(r)+bk<0) \cap (ak-b \ln(r)=0).
\end{equation}
In the case of real $z$, $b=0$, so that the last terms reduces to $k=0$ which implies $n=0$ and $a \ln(r)<0$.  The cut domain is shown in Figure \ref{figure:figure16} as the blue region and the cut trace is shown as the yellow line.

The default secondary branch-cut is thus the line segment  (0,1). 

\subsubsection{Analyzing the $\pLog$ stack}
\begin{figure}
	\centering
			\includegraphics[scale=.5]{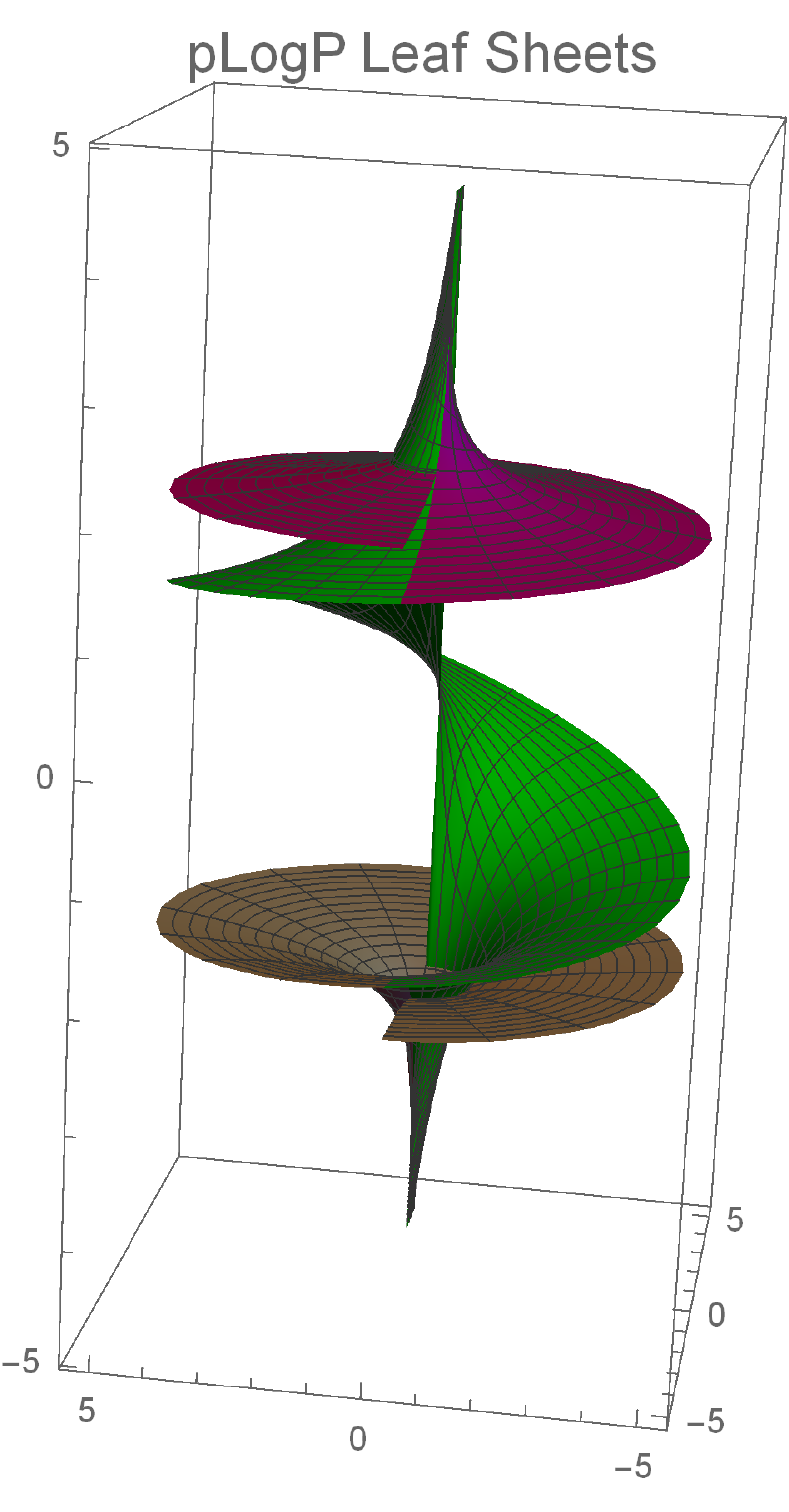}
		\caption{$\pLog$ leaf sheets $\{0,0\},\{-1,0\},\{1,0\}$}
   \label{figure:figure17}
\end{figure}

A section of the $\pLog$ stack for this example is shown in Figure \ref{figure:figure17}.  The primary sheet $\{0,0\}$ is green.  Two secondary sheets of $\{0,0\}$ are also shown:  $\{-1,0\}$ is nickel-color, and $\{1,0\}$ is purple.  These are examples of the default analytic surfaces iterated over by the fixed point iterators in this paper. The combined branch-cut of each sheet is $(-\infty,1)$. This is a Type I log stack.  
And one might suppose the discontinuity over the branch-cut may make the iterator unstable or cyclic.  This is certainly a possibility and has been observed in practice.  However, this can be mitigated by splicing sections of sheets together so that the maximum contiguous surface is presented to the iterator and also by placing the iteration seed close to a root.  However, sheet-splicing was not necessary for the cases tested in this region so the default \texttt{pLog} function was used and simply renamed to \texttt{plog1AN} to easily identify the region begin iterated.
\subsubsection{Creating the auxiliary functions and generating a base contour plot}
The branch and leaf functions are easily computed:

\begin{equation}
\begin{aligned}
\branchF(v,\psi)&=\frac{\Log[\frac{\psi\sec(v\Log(2))}{\Log(2)}]}{\Log(2)},\\
\leafF(v,\phi)&=\frac{\Log[\frac{\phi\csc(v\Log(2))}{\Log(2)}]}{\Log(2)}.
\end{aligned}
\end{equation}

The median of each branch domain is computed using the techniques in the previous section.  $\branchF$ establishes the branch numbers relative to the asymptotes of $\sec(v\Log(2))$. 

\subsubsection{Constructing the Iterator and Computing the roots}
Given
\begin{equation}
T_2(w,n,m;z)=w-\pLog(w,n,m;z)
\end{equation}
with 
\begin{equation}
\pLog(w,n,m)=\frac{1}{\Log(z)}\bigg\{\Log\bigg[\frac{1}{\Log(z)}\left(\Log(w)+2n\pi i\right)\bigg]+2m\pi i\bigg\};\quad (n,m)\in\mathbb{Z},\\
\end{equation}
define:
\begin{equation}
\frac{d}{dw}\pLog(w,n,m)=\frac{1}{w\Log(z)\left(\Log(w)+2n\pi i\right)}.
\end{equation}
A Newton iteration of $T_2$ would take the form:
\begin{equation}
w_{k+1}=w_{k}-\frac{w_{k}-\pLog\left(w_{k},n,m\right)}{1-\frac{1}{w_{k}\Log(z)\left(\Log(w_{k})+2n\pi i\right)}}.
\label{eqn:eqn500}
\end{equation}
Note the branch parameters  were computed symbolically with arbitrary precision.  Therefore the seeds have arbitrary precision as well.  The general procedure for iterating the roots using (\ref{eqn:eqn500}) is the following:
\begin{enumerate}
\item
Set a working precision for the seed, and a desired accuracy for the root.  The seed precision is usually set between 30 and 100 digits.
\item
Begin iterating:  Iterate until either a desire accuracy is achieved or a maximum number of iterations is reached.  
\item  
Due to the close proximity of sequential roots on a branch lobe, if sequential roots are computed, the seed for the next root is set to the previously calculated root.  However, the previous root is first rationalized to a rational number prior to setting it's precision to the starting seed precision.
\end{enumerate}

\subsubsection{Test Results}
Table \ref{table:table2} lists roots $\{1,0\}$ to $\{10,0\}$ for $z=2$  computed with $\texttt{pLog1AN}$ using a working precision of $40$ and accuracy of $30$.  The roots are accurate to $27$ decimal digits.  Column "`I"' is the number of iterations used to achieve this accuracy.  Computation of the first $100,000$ roots with an accuracy of $47$ digits on a $2.2$ GHz machine took $55$ seconds. 

\begin{table}[!ht]
\small
\centering
\begin{tabular}{|c|c|c|}
\hline
Root & I&Value \\
\hline
1	&5&	$3.32659640338718651629111245784+1.98345940154898674596183895638i$\\
2	&4&	$4.24218210955016856138327332701+2.09484399760721115618994484258 i$\\
3 &4&	$4.80231963088007566818857548380+2.14221283941333326231852002691 i$\\
4	&4&	$5.20604757038160215682362465030+2.16856760479743763361704916320 i$\\
5	&4&	$5.52165579954147892613237150892+2.18542615221552883665073000990 i$\\
6	&4&	$5.78070954903593282245775341453+2.19717276951965471971748103016 i$\\
7	&3&	$6.00038886534664208553763234121+2.20584495530286752551209605520 i$\\
8	&3&	$6.19107830256296787260245128971+2.21252055199734402929507988891 i$\\
9	&3&	$6.35953349646166768755094539015+2.21782426708596179600889292447 i$\\
10	&3&	$6.51039555600714279854298752089+2.22214366914367902544431058445 i$\\
\hline
\end{tabular}
\caption{Roots $\{1,0\}$ to $\{10,0\}$}
 \label{table:table2}
\end{table}

\begin{figure}[!ht]
	\centering
			\includegraphics[scale=.5]{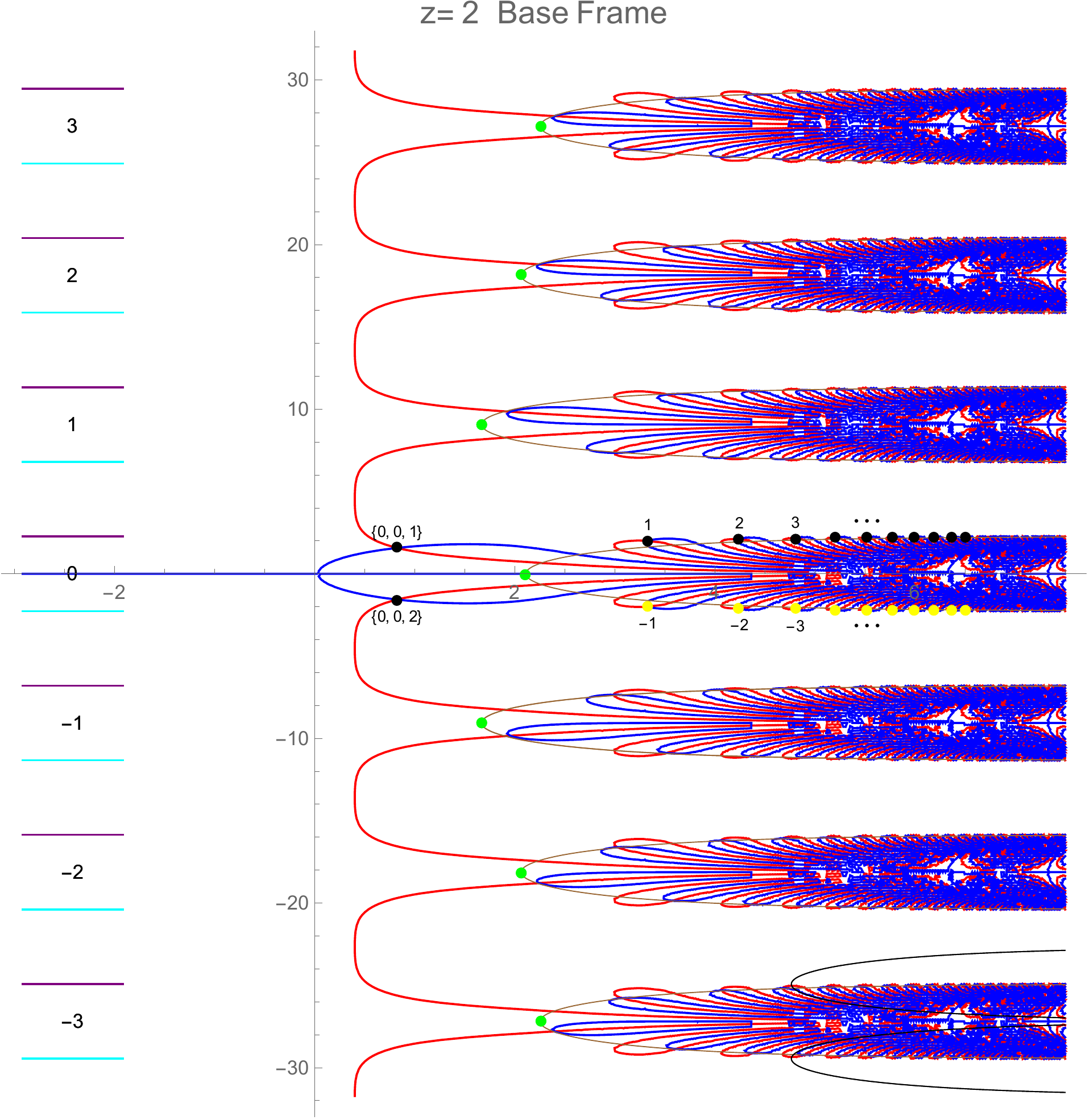}
		\caption{Contour Diagram for Region 1A}
   \label{figure:figure20}
\end{figure}

Figure \ref{figure:figure20} is a contour plot for this example.  Most roots are arranged on the leaf bulbs.  Branches near the origin may have isolated roots with \quotes{degenrate} leaf lobes on the real axis as shown in Figure \ref{figure:figure15}.  Seven branch lobes are shown in the diagram. The roots in Table \ref{table:table2} are shown as black points.   Twenty roots are shown on branch $0$ as the black and yellow points, several of which are labeled. Example plots of the auxiliary function $\branchF$ are the brown traces encompassing the branch bulbs.  The diagram shows two black traces of $\leafF$ on branch $-3$.  The constant $\phi$ in $\leafF$ can be adjusted so that the trace passes close to a root.  Since  $\leafF$ is an analytic expression, a seed very close to a root could be computed using the derivative of $\leafF$.  However, this was not needed in this study.

Contour plots show only a small number out of an infinite number of $T_2$ branch and leaf bulbs.  Although the distance between branches is a constant, the distance between successive leaf lobes decrease.   For example, root $\{10^{12},0\}$ of this example has real part near $43$.

Seeds used to seed the $\pLog$ iterators are shown in a contour plot as green points close to the head of each branch bulb.   The seeds are computed using $\branchF$ as described in the previous section.   The branch numbers are shown between purple and cyan lines delimiting each branch. 
 
The orientation of the branching in Figure \ref{figure:figure20} is a reflection that it is in normal form.   Roots lying on the same branch are grouped together and identified by $\{n,m,p\}$ for constant $m$ with $p$ omitted if the root is not a multiple root.

Note roots $\{0,0,1\}$ and $\{0,0,2\}$ in Figure \ref{figure:figure20}.  This is an example of a $\pLog$ sheet having two roots:  Iteration of $\text{pLog1AN}$  over a $20\times20$ region in the $w$-plane centered at the origin for root $\{0,0\}$ produces two basins:  The upper half-plane attracts root $\{0,0,1\}$ and the lower half-plane, $\{0,0,2\}$. 

Table \ref{table:table3} lists roots $\{10^{12},10^{12}\}$ to $\{10^{12}+10,10^{12}\}$ computed with an accuracy of $37$ digits.  The working precision was set to 50 digits.  Note only one iteration was needed to achieve $37$ digits of accuracy.   Iterating the values once more raised the accuracy to $85$ digits.

\begin{table}[!ht]
\small
\centering
\begin{tabular}{|c|c|c|}
\hline
Root & I& Value \\
\hline
1	&1&	43.04339964106736795553573548+9.06472028365665379932627263776832711482*$10^{12}i$\\
2	&1&	43.04339964106881065057662192+9.06472028365665379932627263776832712167*$10^{12}i$\\
3	&1&	43.04339964107025334561750692+9.06472028365665379932627263776832712853*$10^{12}i$\\
4	&1&	43.04339964107169604065839047+9.06472028365665379932627263776832713538*$10^{12}i$\\
5	&1&	43.04339964107313873569927258+9.06472028365665379932627263776832714223*$10^{12}i$\\
6	&1&	43.04339964107458143074015325+9.06472028365665379932627263776832714908*$10^{12}i$\\
7	&1&	43.04339964107602412578103247+9.06472028365665379932627263776832715593*$10^{12}i$\\
8	&1&	43.04339964107746682082191025+9.06472028365665379932627263776832716278*$10^{12}i$\\
9	&1&	43.04339964107890951586278659+9.06472028365665379932627263776832716963*$10^{12}i$\\
10	&1&	43.04339964108035221090366149+9.06472028365665379932627263776832717648*$10^{12}i$\\
\hline
\end{tabular}
\caption{First 10 positive roots over one-trillion on branch one-trillion, $z=2$}
 \label{table:table3}
\end{table}
%
%
\subsection{Region 1B}
Region 1B has a root of multiplicity $2$ at $w=e$ so an over-relaxation factor of $2$ is included in the iteration method for this root: 
\begin{equation}
w_{k+1}=w_{k}-2\left(\frac{w_{k}-\pLog\left(w_{k},n,m\right)}{1-\frac{1}{w_{k}\Log(z)\left(\Log(w_{k})+2n\pi i\right)}}\right).
\end{equation}
Since this region has the same branch-cuts as Region 1A, $\pLog$ was used in all test cases and produced roots to the desired accuracy.  The iterator name is changed to $\texttt{pLog1BN}$.  
%
%
\subsection{Region 1C}
This region has the same branch-cuts as 1A and 1B and can be iterated with $\pLog$.  However, the iterator splits into two basins for root $\{0,0\}$ corresponding to the two real roots $\{0,0,1\}$ and $\{0,0,2\}$ in Figure \ref{figure:figure15}. 
%
%
\subsection{Region 1D}

\begin{figure}[!ht]
	\centering
			\includegraphics[scale=.5]{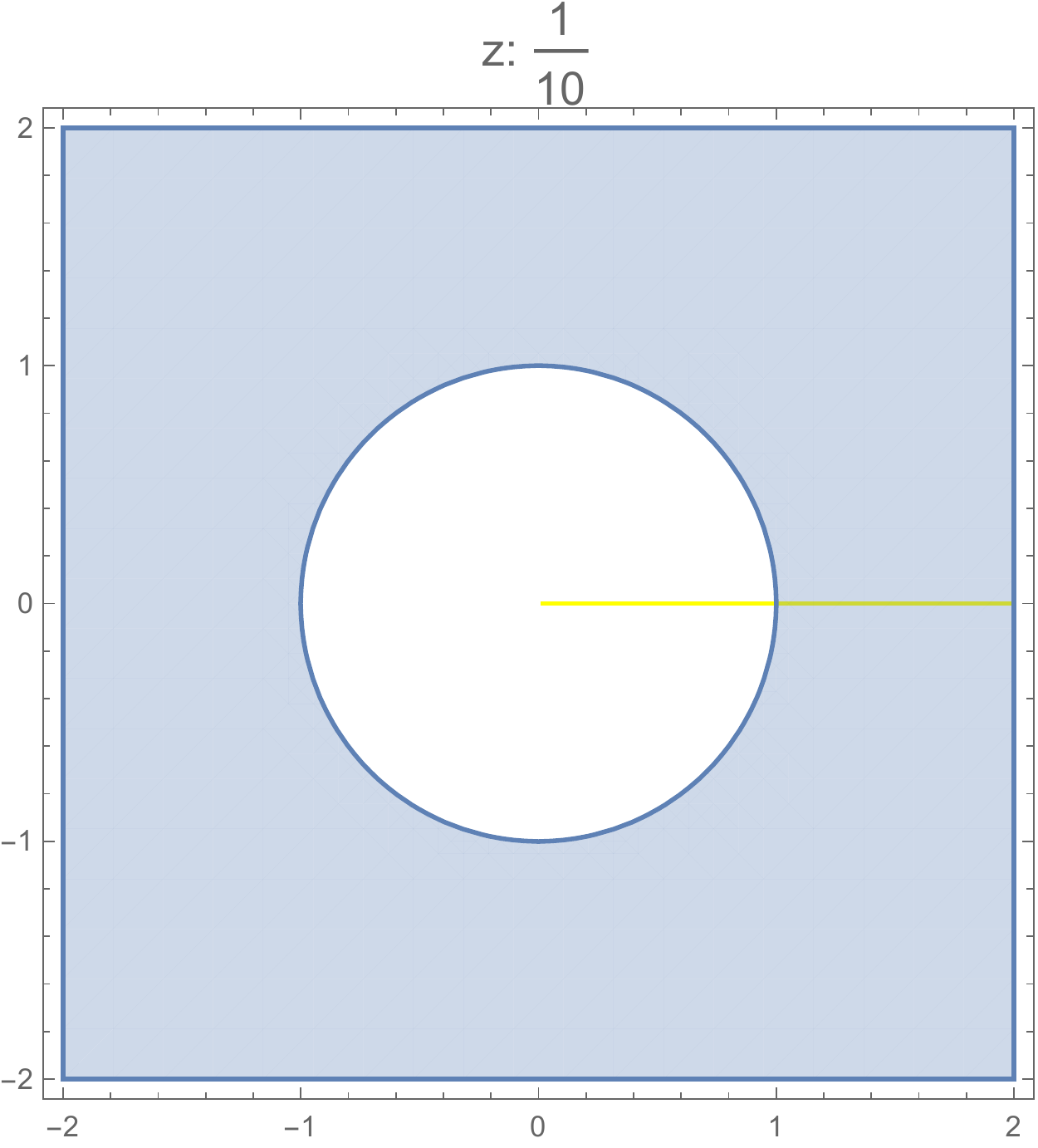}
		\caption{Contour Diagram for Region 1D}
   \label{figure:figure24}
\end{figure}

Figure \ref{figure:figure24} is a region plot of the secondary branch- cut.   Note this time, the intersection of the cut domain and trace is now the line $(1,\infty)$. 
\begin{figure}
	\centering
			\includegraphics[scale=.5]{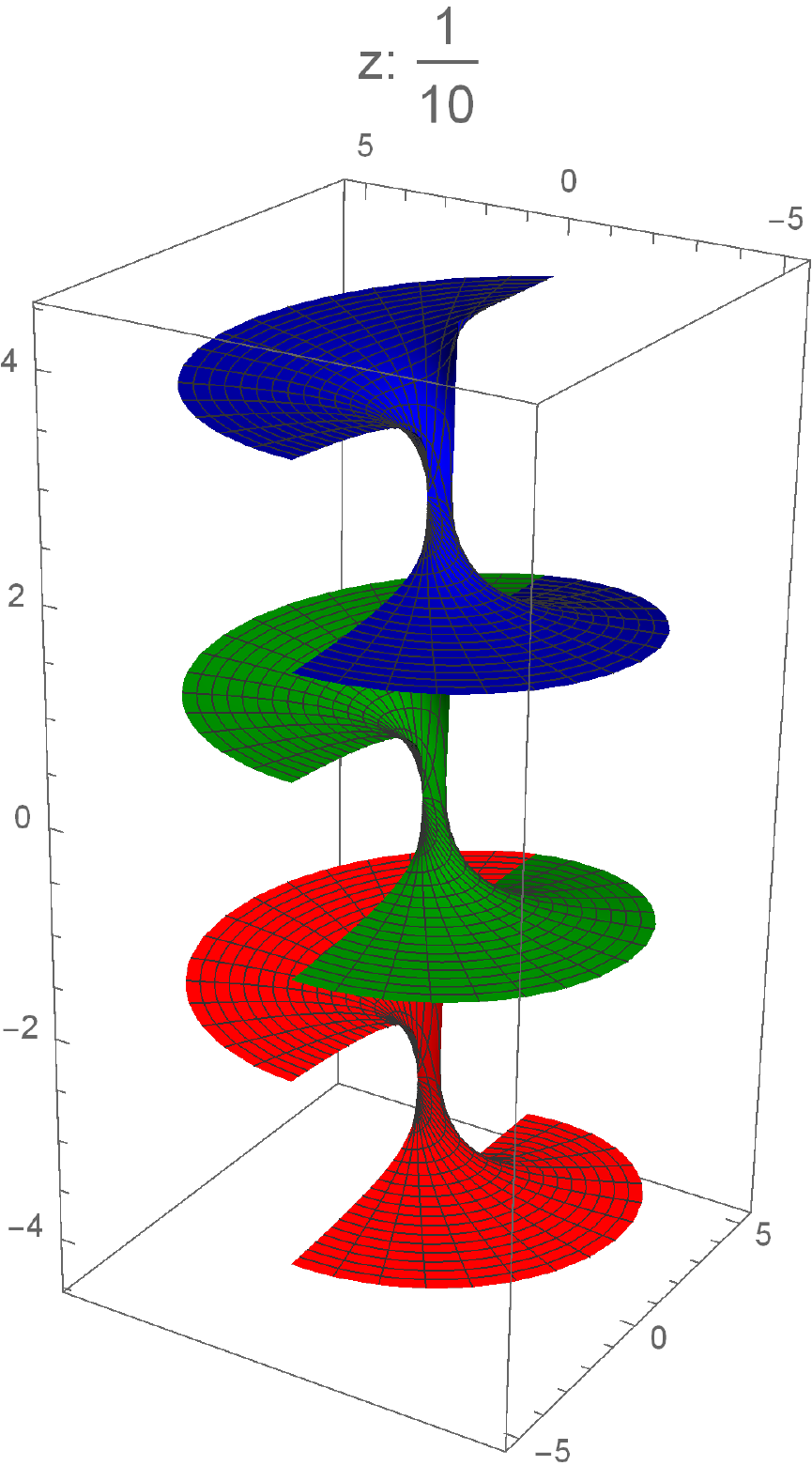}
		\caption{Type II log stack for Region 1D}
   \label{figure:figure25}
\end{figure} 

The secondary cut results in a Type II log stack shown in Figure \ref{figure:figure25}.  The primary $\pLog$ sheet $\{0,0\}$ is green, $\{0,-1\}$ is blue and $\{0,1\}$ is red.  In order to align root indexes with the contour branches, the $\pLog$ stack is converted to a  Type I stack with branch-cut   $(-\infty,1)$.   The code to do this is shown in Listing \ref{code1}. 

\begin{center}
\begin{minipage}{0.75\linewidth}
\begin{lstlisting}[language=Mathematica,
frame=single,
caption=Mathematica code for the pLog1DN iterator,
label=code1]
pLog1DN[w_, n_, m_] := Module[{logVal},
   If[n == 0,
	(*
	 for {0,0} root, use upper half section of branch 0 with lower half of branch -1
	*)
	  If[Arg[w] >= 0,
      logVal=pLog[w, n, m];
      ,
      logVal=pLog[w, n, m - 1];
      ];
     ,
		(* 
		use branch m-1 for negative roots and branch m for positive roots 
		*)
     If[n < 0,
      logVal=pLog[w, n, m - 1];
      ,
      logVal=pLog[w, n, m];
      ];
     ];
    logVal
   ];
\end{lstlisting}
\end{minipage}
\end{center}

\begin{figure}[!ht]
	\centering
			\includegraphics[scale=.5]{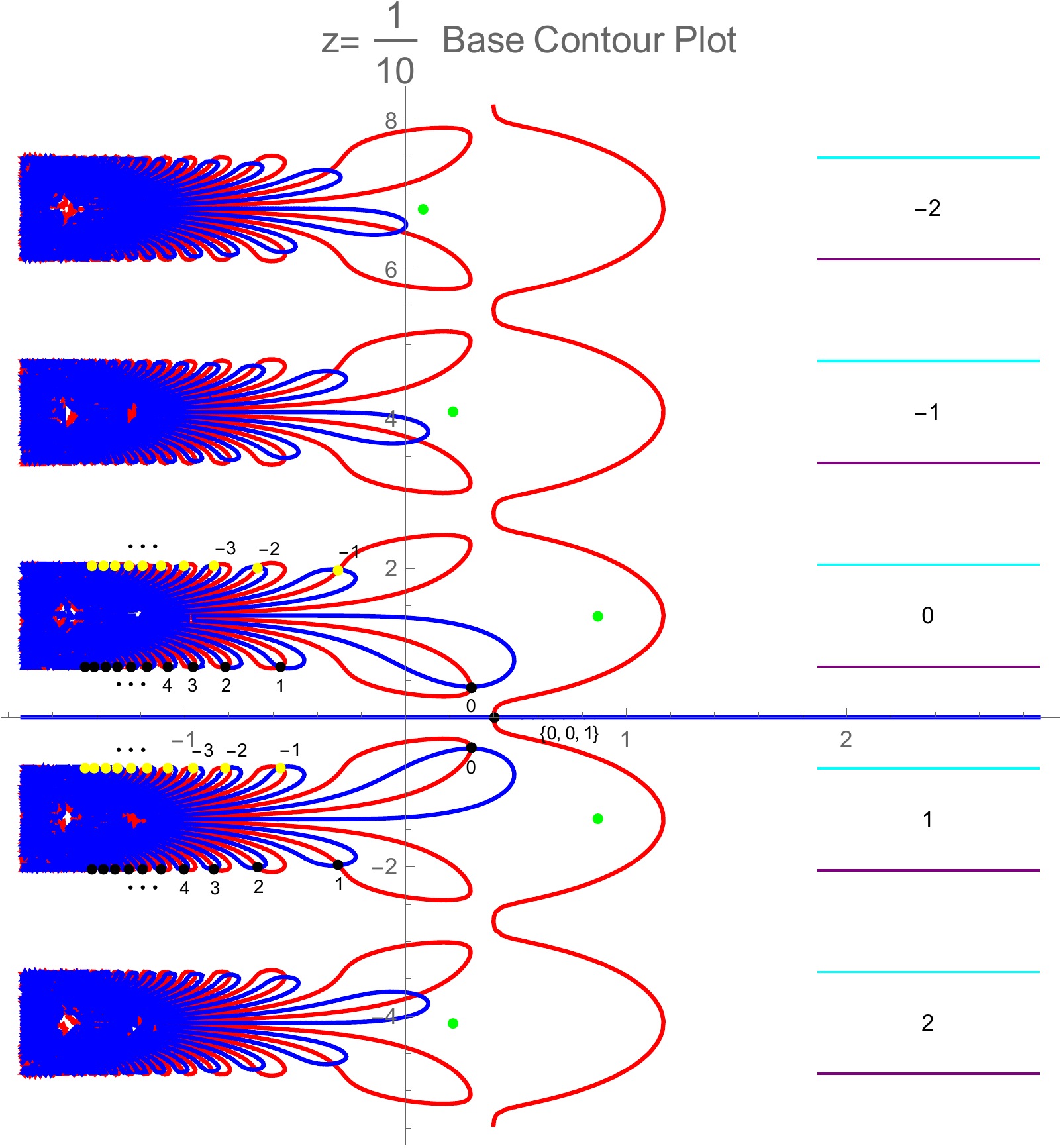}
		\caption{Contour Diagram for Region 1D}
   \label{figure:figure26}
\end{figure}

\texttt{pLog1DN} stitches the upper half section of branch $\{0,m\}$ with the lower section of $\{0,m-1\}$ while iterating over $\{n,m\}$ for positive roots and $\{n,m-1\}$ for negative roots.  This code effectively converts the default log stack to a Type I stack.  For the single real root, the default $\pLog$ iterator is used to iterate over the default (green) $\{0,0\}$ sheet.  This is done because the real root is at the branch-cut for the  Type I log sheets of \texttt{pLog1DN}  whereas it is in the center of the Type II $\pLog$ sheet surrounded by an analytic domain and so less affected by a branch-cut.   Figure \ref{figure:figure26} is an example contour diagram of this region.  Note the branch lobes are horizontal and opening towards the left.  In order to generate the normal frame of this region as well as Regions 1E and 1F, a rotation transformation with $\beta$ set to $-\pi$ is applied.  

\begin{figure}[!ht]
	\centering
			\includegraphics[scale=.5]{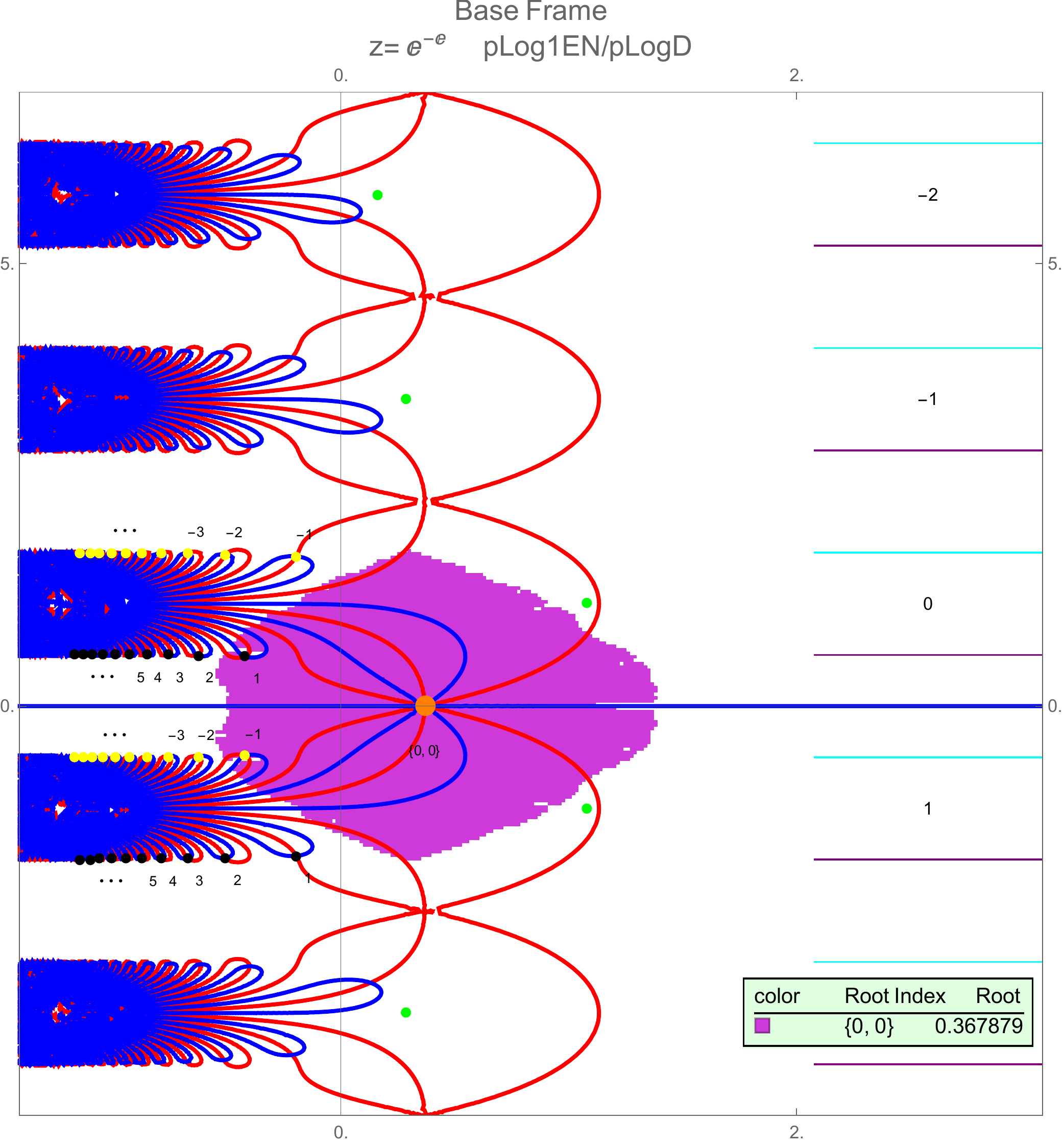}
		\caption{Region 1E contour diagram with basin for root $\{0,0\}$}
   \label{figure:figure27}
\end{figure}
%
%
\subsection{Region 1E}

The Type I iterator of Region 1D can be used for Region 1E and named $\texttt{pLog1EN}$.  The real root at $1/e$ is assigned to branch $\{0,0\}$ because the root is in the region of definition of \texttt{pLog1EN} and so labeled $\{0,0,1\}$. 

Figure \ref{figure:figure27} is a contour diagram superimposed on the basin diagram for root $\{0,0\}$.  The orange point at $1/e$ is the real root $\{0,0\}$ and has multiplicity $3$ so a relaxation factor of $3$ is used in the iterator for this root.  The purple region in the figure is the basin of convergence of this root.  The iterator diverges in the surrounding white area.  The title of the diagram \quotes{pLog1EN/pLogD} means $\pLog$ was used to iterate root $\{0,0\}$ with the accompanying basin diagram, and $\texttt{pLog1EN}$ was used to iterate the others.  

%
%
\subsection{Region 1F}

\begin{figure}
	\centering
			\includegraphics[scale=.5]{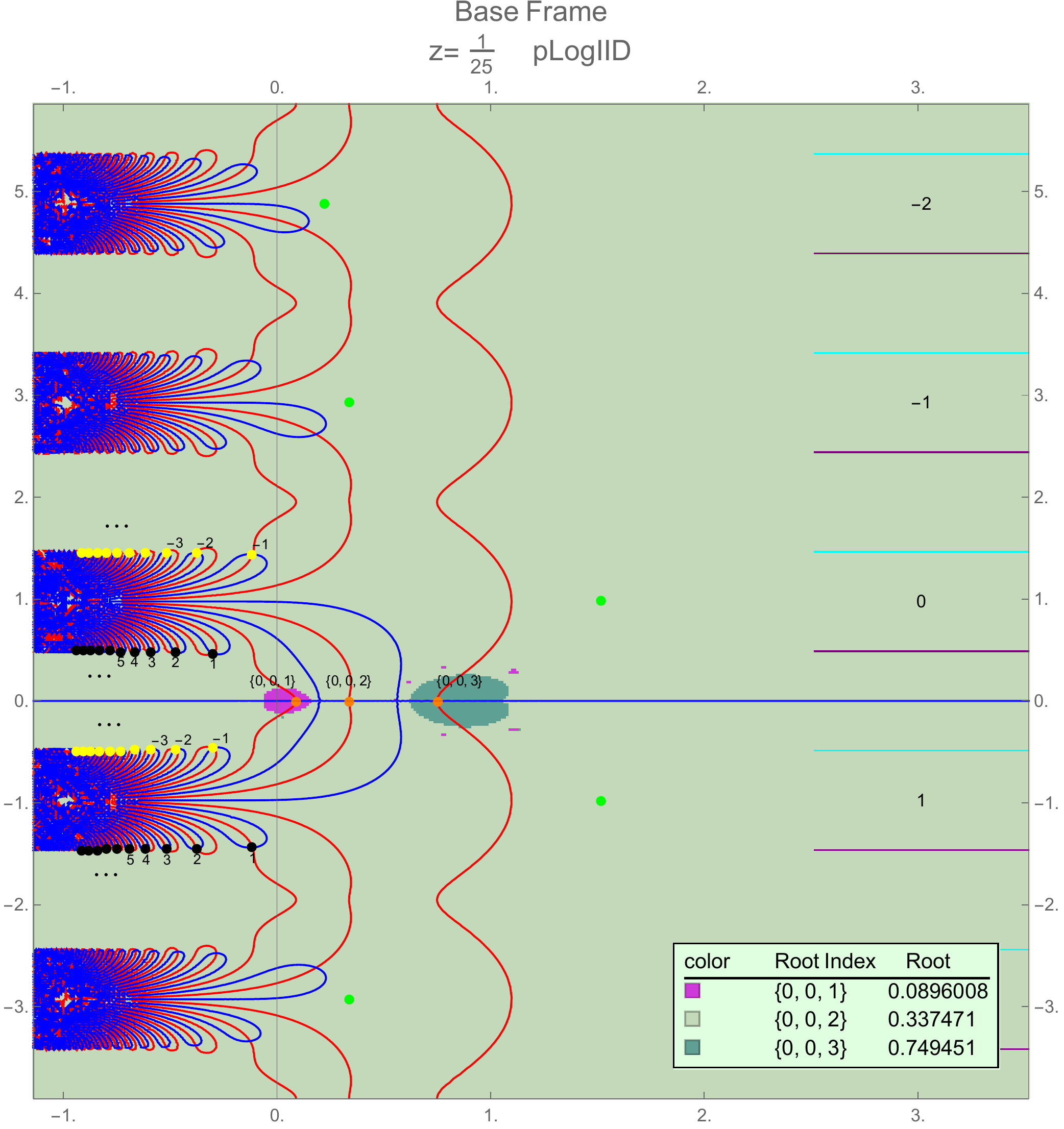}
		\caption{Contour diagram and basins for Region 1F}
   \label{figure:figure28}
\end{figure}

The branch-cuts for Region 1F are the same as Region 1E.  All roots except the three real roots can be iterated with $\texttt{pLog1EN}$  renamed $\texttt{pLog1FN}$.  The real roots are iterated with $\pLog$ which has the basin diagram shown in the background of Figure \ref{figure:figure28}. 

As $z$ approaches the origin, the basins for roots $\{0,0,0\},\{0,0,1\}$ and $\{0,0,2\}$ in Region 1F grow smaller making it more difficult to locate a seed for these roots.  However, one could incrementally approach these basins by starting the iteration with a relatively large value of z, using the associated roots as seeds for smaller values of z until the required values are reached.
%
%
\section{Region 2A}

Region 2 is the unit circle: 2A is the upper half, and 2B is the lower half.  For $z=e^{i\theta}$, the branching lobes open downward when $\theta>0$ and upward when $\theta<0$.  Therefore, the rotation angle $\beta$, is  $\pm \pi/2$. Recall the definitions
$$
\begin{array}{cc}
\text{cut Domain:}& \ln(r)+bk<0, \\
\text{cut Trace: }& ak-b\ln(r)=0
\end{array}
$$

with the secondary branch-cut the intersection $(a\ln(r)+bk<0)\cap (ak-b \ln(r)=0)$.  Two cases are considered:
\begin{enumerate}
\item
$n=0$:  In the case of branch $\{0,0\}$, $b\theta<0$ for the first requirement.  This is the lower half-plane.  In the second term, $b\ln(r)=0$.  This is the unit circle so thus the intersection of these is the unit half-circle in the lower plane.   Figure \ref{figure:figure29a} is a plot of $\pLog$ branches with $\{0,0\}$ green, $\{0,1\}$ red,  and $\{0,-1\}$ blue.  Note the default cut $(-\infty,0)$ and the half-circle cut in the lower half-plane.

The next objective is to stitch together components of the branches  creating  a derived \texttt{pLogR2N} stack with a single branch-cut.   Consider first the green $\{0,0\}$ branch sheet in Figure \ref{figure:figure29a}.  There are two ways to stitch this sheet with the others to create a contiguous sheet with a single branch-cut $(-\infty,1)$:
\begin{enumerate}
\item
Combine the green sheet outside the lower half-circle with the blue sheet inside the half circle,
\item
Combine  the green sheet inside the lower half-circle with the red sheet outside the lower half-circle,
\end{enumerate}
  
\item
$n\ne 0$:  We then have $b(\theta+2n \pi)<0$ for the first requirement.  For $n>0$ this is never true so that the positive branches have only the default cut.  However, if $n<0$, then this is always true. In the second term, $b \ln(r)=0$, which again gives us the unit circle so that the branch-cuts for the negative sheets is the unit circle with the primary branch-cut.  The log stack for this case is shown in Figure \ref{figure:figure29b} with $\{-1,0\}$ brown, $\{-1,1\}$ yellow, and $\{-1,-1\}$ purple.

\begin{figure}[!ht]
     \centering
     \begin{subfigure}{0.4\textwidth}
         \centering
         \includegraphics[width=0.75\textwidth]{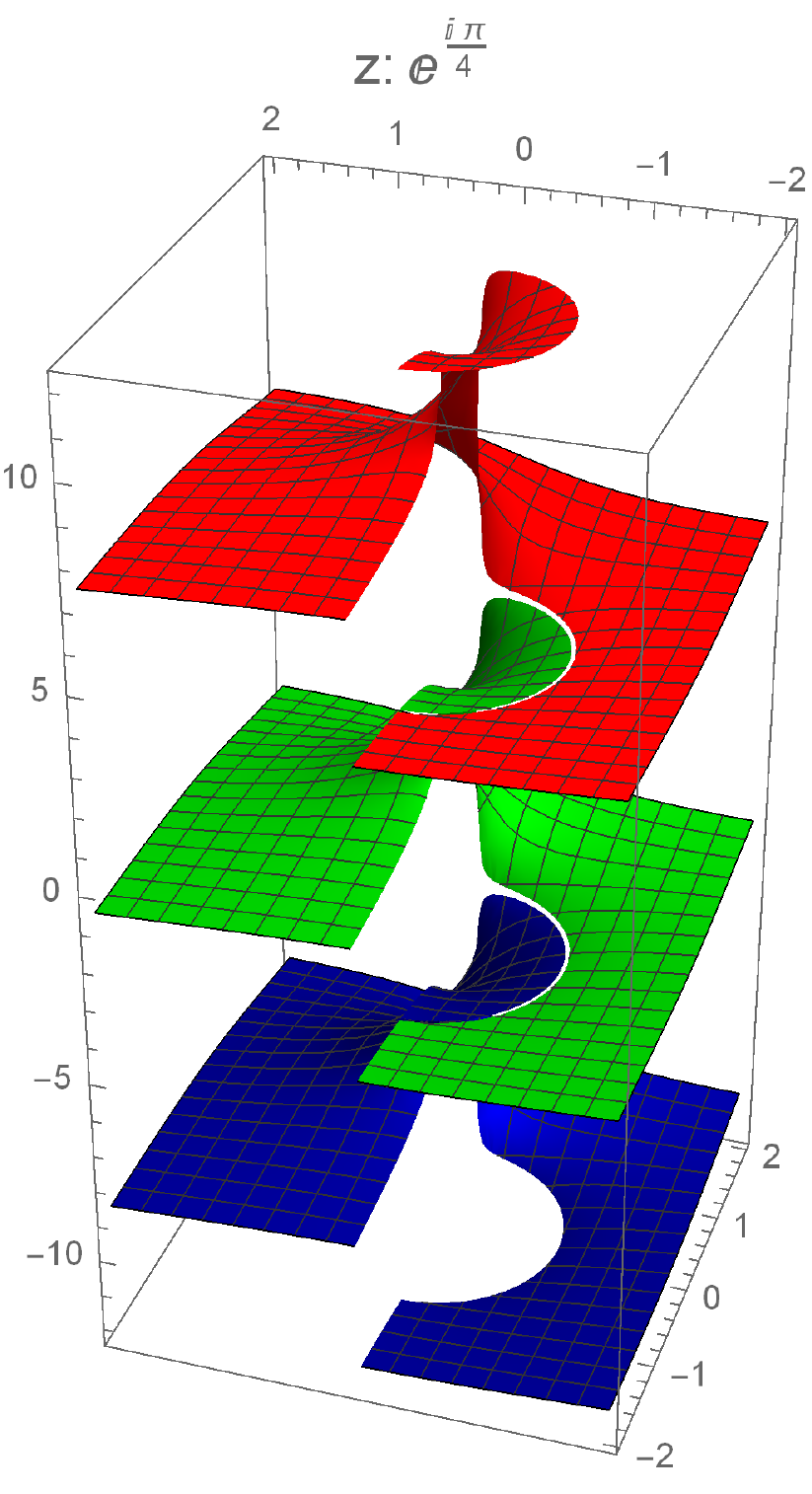}
         \caption{Case 1 sheets $\{0,0\},\{0,1\},\{0,-1\}$}
         \label{figure:figure29a}
     \end{subfigure}
     \hfill
     \begin{subfigure}{0.4\textwidth}
         \centering
         \includegraphics[width=0.75\textwidth]{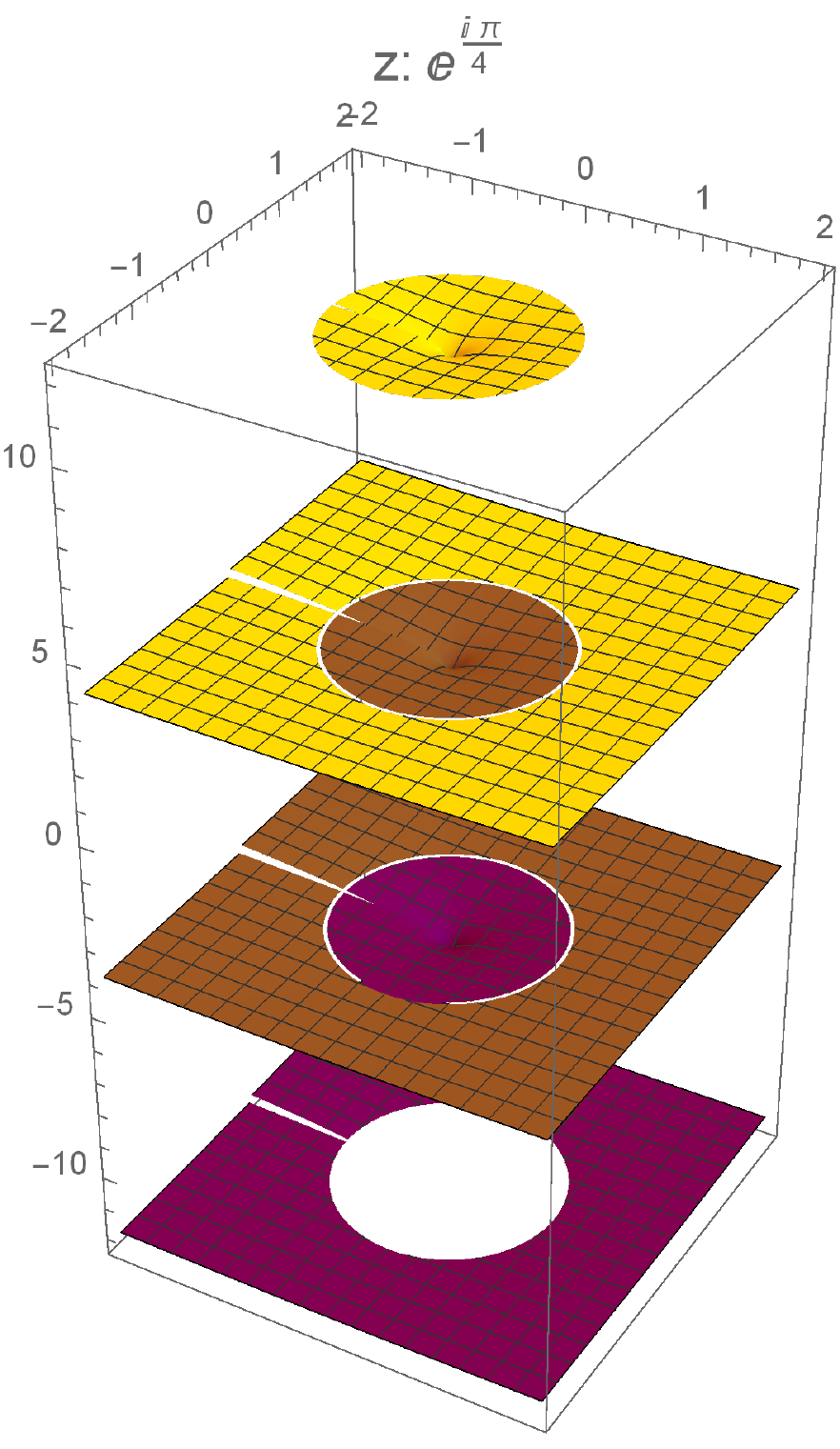}
         \caption{Case 2 sheets $\{-1,0\},\{-1,1\},\{-1,-1\}$ }
         \label{figure:figure29b}
     \end{subfigure}
     \hfill
     \caption{Region 2A branch-cuts}
        \label{figure:figure29}
\end{figure}

Now consider the brown $\{-1,0\}$ sheet in Figure \ref{figure:figure29b}.  There are two ways to combine this sheet with the others to create a contiguous sheet with branch-cut $(-\infty,0)$:
\begin{enumerate}
\item
Combine the brown sheet outside the unit circle with the purple sheet inside the unit circle,
\item
Combine the brown sheet inside the unit circle with the yellow sheet outside the unit circle.
\end{enumerate}

In order to align the roots with the contour plot, 1a and 2a are combined to create derived stack \texttt{pLog2AN}.  The code to construct this stack is shown in Listing \ref{code2}.
\end{enumerate}
\begin{figure}[!ht]
	\centering
			\includegraphics[scale=.65]{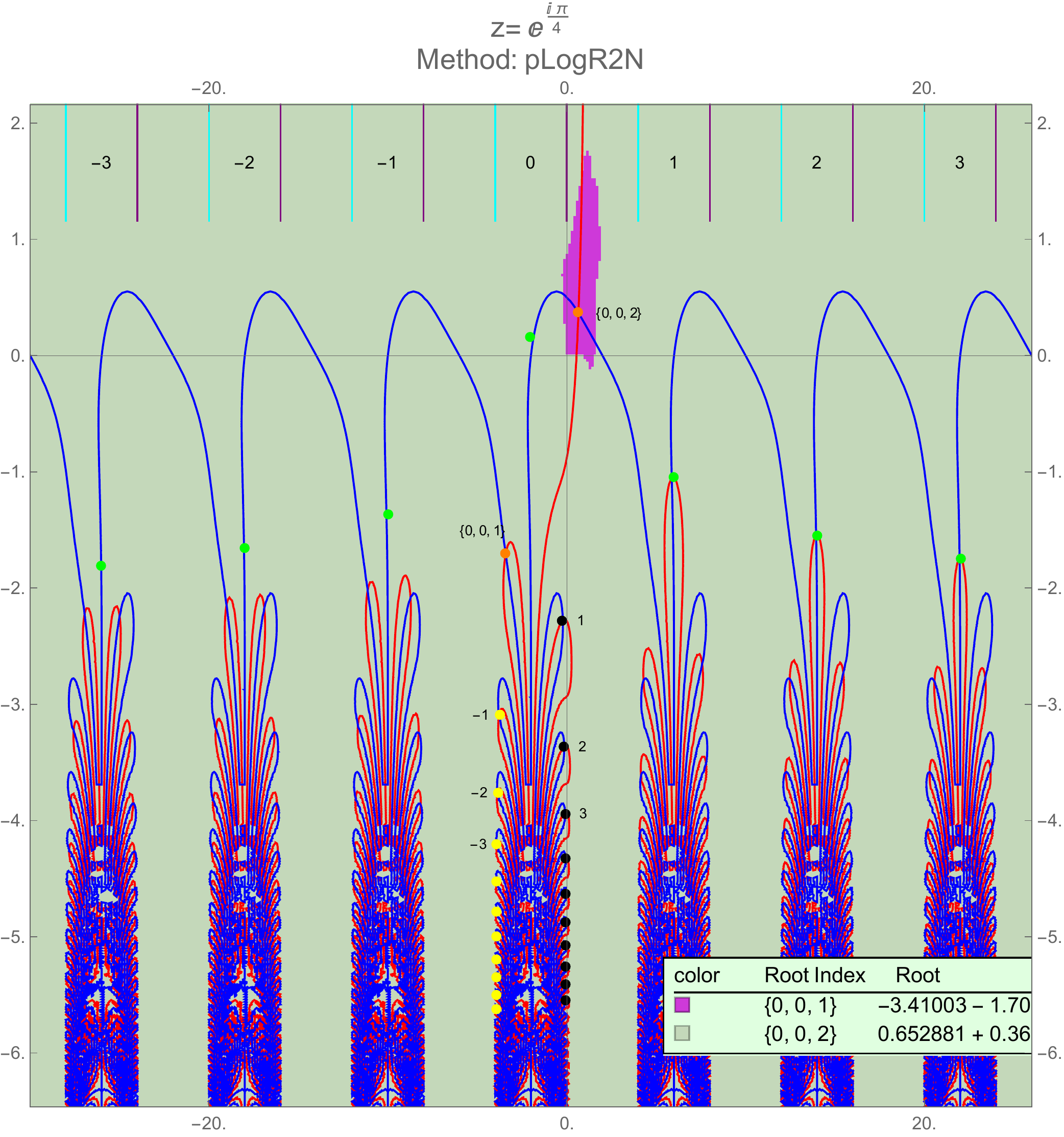}
		\caption{Region 2A contour plot with basins for $\{0,0\}$}
   \label{figure:figure33}
\end{figure}
\begin{center}
\begin{minipage}{0.75\linewidth}
\begin{lstlisting}[language=Mathematica,
frame=single,
caption=Mathematica code for the pLog2AN iterator,
label=code2]
pLog2AN[w_, n_, m_] := Module[{pVal},
   If[n == 0,
    If[Arg[w] >= 0,
      pVal = pLog[w, n, m];
      ,
      If[Abs[w] > 1,
        pVal = pLog[w, n, m];
        ,
        pVal = pLog[w, n, m - 1];
        ];
      ];
    ,
    If[n > 0,
      pVal = pLog[w, n, m];
      ,
      If[Abs[w] >= 1,
        pVal = pLog[w, n, m];
        ,
        pVal = pLog[w, n, m - 1];
        ];
      ];
    ];
   pVal
   ];
\end{lstlisting}
\end{minipage}
\end{center}
Figure \ref{figure:figure33} is a contour plot with $10$ negative and $10$ positive roots of branch $0$ as the black and yellow points.  A basin diagram for root $\{0,0\}$ is in the background.  Roots $\{0,0,1\}$ and $\{0,0,2\}$ corresponding to the basins are shown as the orange points, and the branch iteration seeds as the green points.   A seed in the light green area will iterate to root $\{0,0,2\}$.  A seed in the small purple region will iterate to root $\{0,0,1\}$.  The roots were computed to $43$ digits of accuracy in under $6$ iterations using \texttt{pLog2AN}. 
%
%
\section{Region 2B}
Region 2B is the lower half unit disc.  The branch-cuts in this region is the inverse of Region 2A:  the intersection of the principal cut with the lower half unit circle.  The root-finding algorithms are equivalent to the Region 2A case except now we consider the half disc in the lower half plane.  The same techniques of Region 2A can be used to construct \texttt{pLog2BN} for this region.
%
%
\section{Region 3A}
\begin{figure}[!ht]
     \centering
     \begin{subfigure}{0.4\textwidth}
         \centering
         \includegraphics[width=0.75\textwidth]{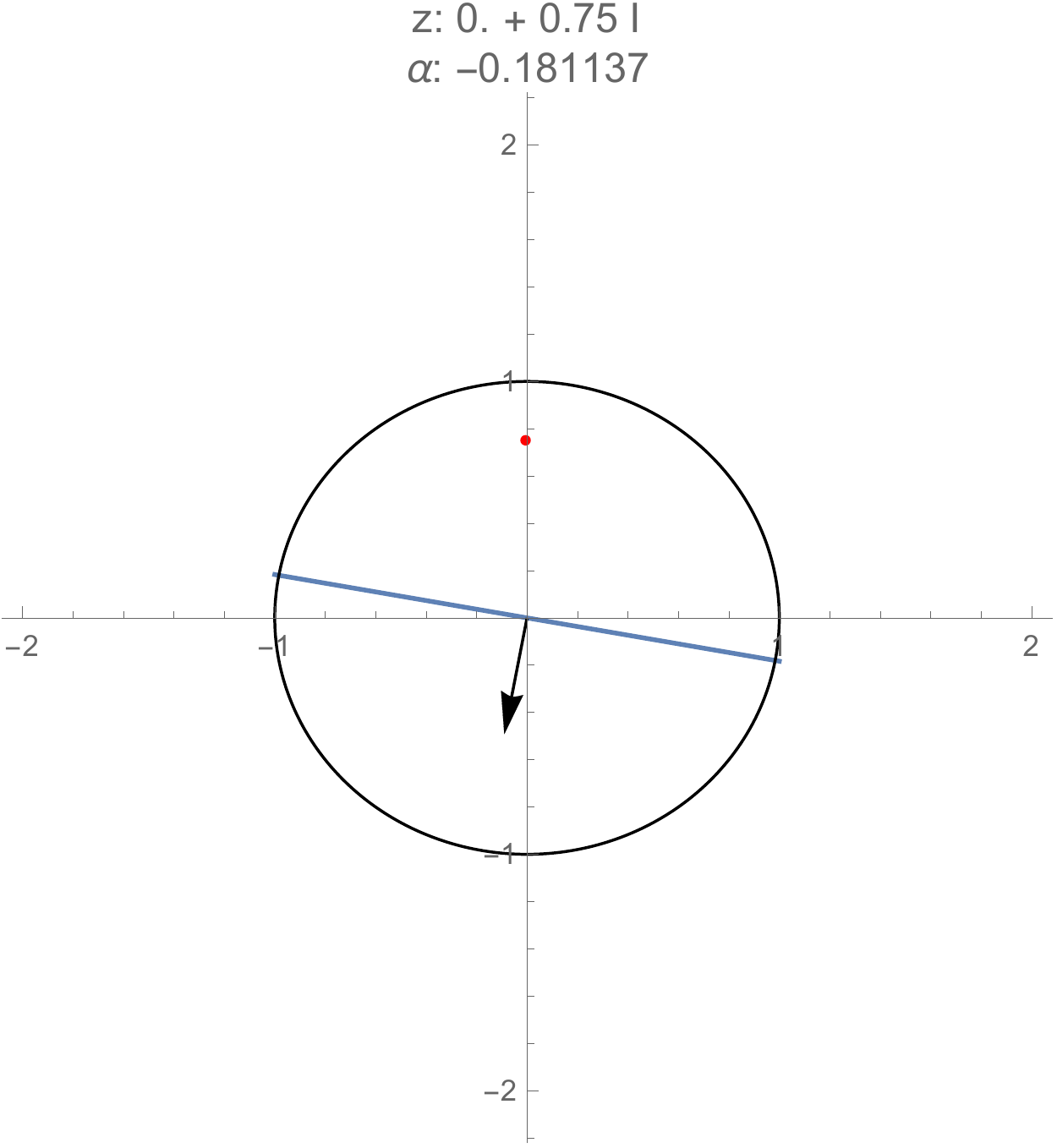}
         \caption{$z=3/4i$}
         \label{figure:figure34b}
     \end{subfigure}
     \hfill
     \begin{subfigure}{0.4\textwidth}
         \centering
         \includegraphics[width=0.75\textwidth]{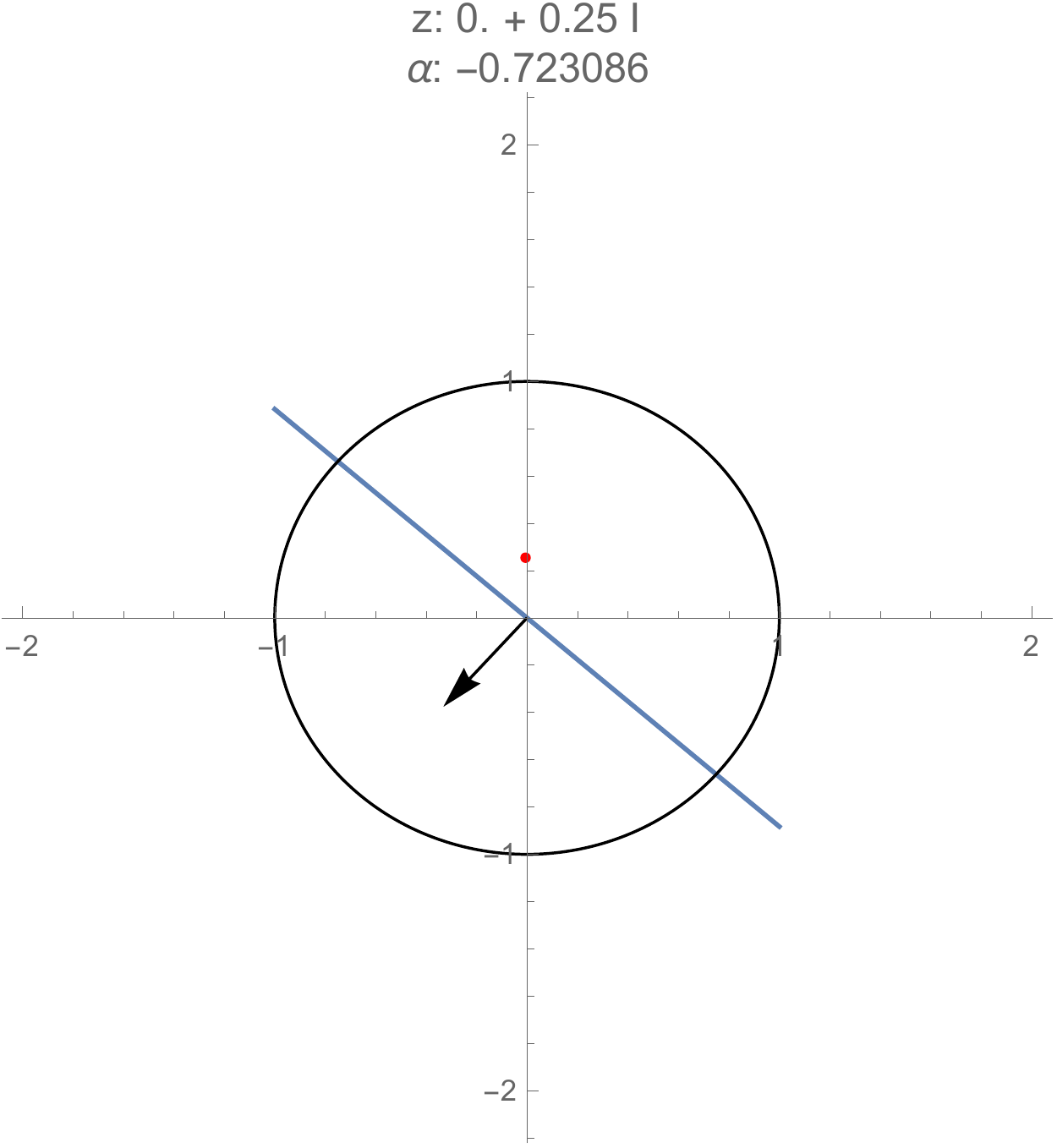}
         \caption{$z=1/4i$}
         \label{figure:figure34a}
     \end{subfigure}
     \hfill
     \caption{Region 3A Inclination Angles}
        \label{figure:figure34}
\end{figure}

The branching and leaf contours of Region 3A open downward with an inclination angle that becomes more negative as z approaches the origin.  Two cases of the inclination axis are shown in Figure \ref{figure:figure34} with the red point indicating the location of $z$.  The branch leaf lobes open in the direction of the arrows.

 In Figure \ref{figure:figure34b}, $z$ is close to the unit circle and therefore close to Region 2A so the inclination angle is close to zero.  As $z$ approaches zero, the inclination axis becomes more negative as shown by the second plot.  The angle $\Arg(z)$ has only a small effect on the inclination angle in this region. 

\begin{figure}
     \centering
     \begin{subfigure}{0.4\textwidth}
         \centering
         \includegraphics[width=0.75\textwidth]{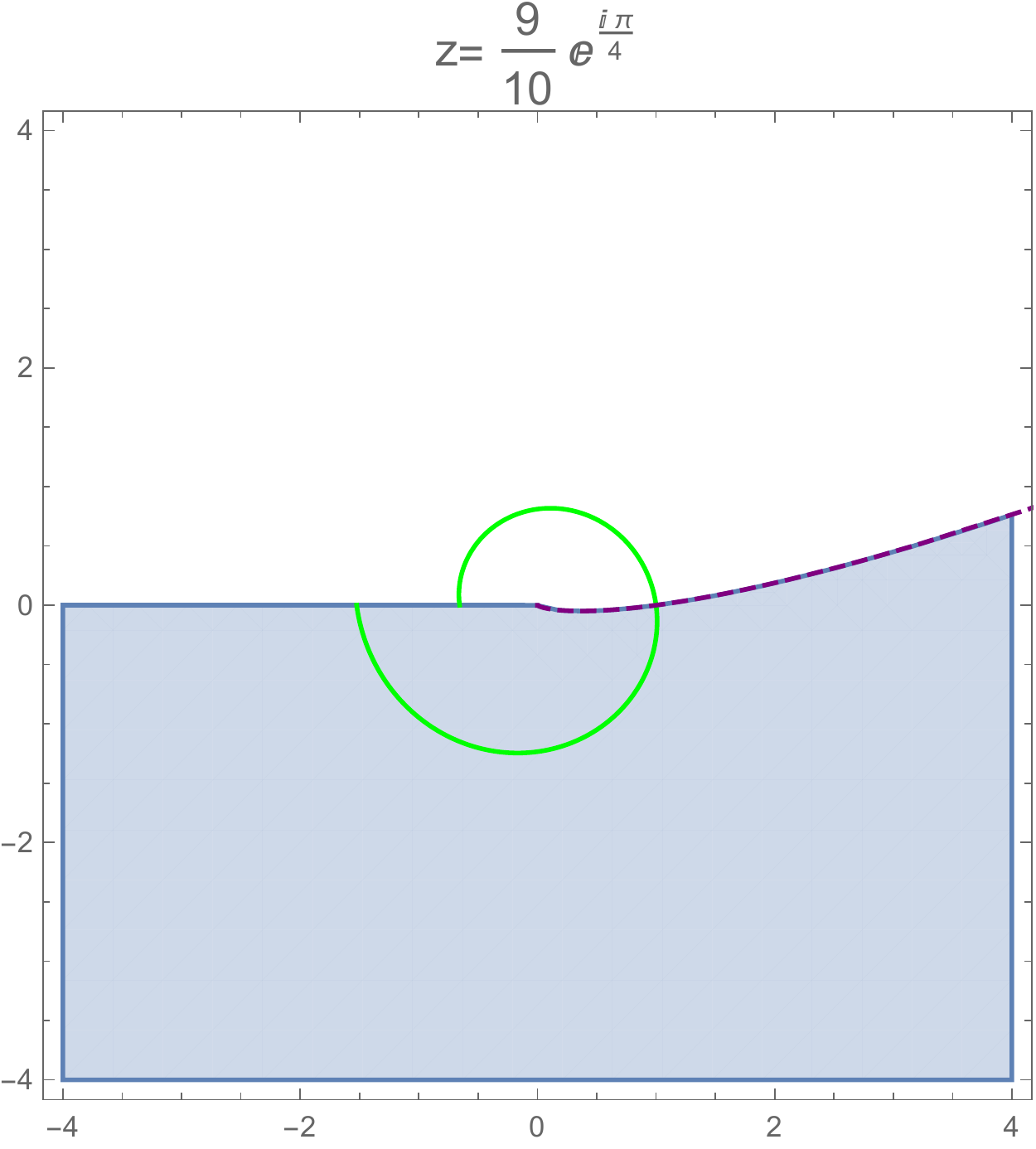}
         \caption{Secondary branch-cut for $\{0,m\}$}
         \label{figure:figure35a}
     \end{subfigure}
     \hfill
     \begin{subfigure}{0.4\textwidth}
         \centering
         \includegraphics[width=0.75\textwidth]{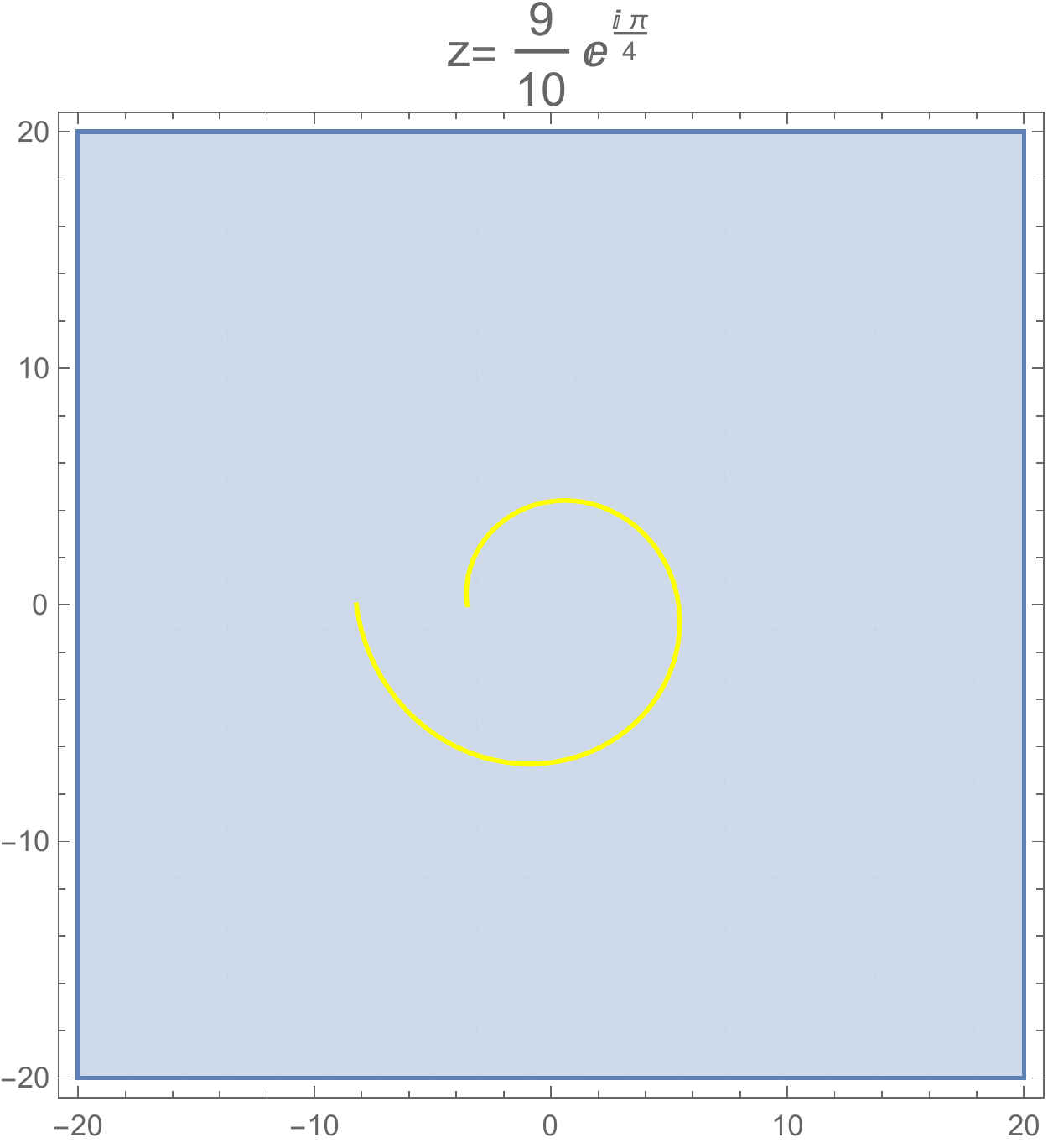}
         \caption{Secondary branch-cut for $\{-1,m\}$}
         \label{figure:figure35b}
     \end{subfigure}
     \hfill
     \caption{Region 3A Secondary Branch-Cuts}
        \label{figure:figure35}
\end{figure}

The secondary branch-cuts however are more complicated in regions 3 and 4. Recall, the secondary branch-cut is the intersection $(a\ln(r)+bk<0)\cap(ak-b \ln(r)=0)$. Re-writing these expressions as
\begin{equation}
\begin{aligned}
\text{cut domain}&:r&>e^{-b(\theta+2 n\pi)}, \\
\text{cut Trace}&:r&=e^{a/b(\theta+2 n\pi)},
\end{aligned}
\end{equation}
show the domain boundary and trace are logarithmic spirals.  First consider the case of branch $\{0,m\}$ shown in Figure \ref{figure:figure35a} with $z=9/10 e^{\pi i/4}$. This plot shows the secondary cut domain in blue and the cut trace in green.  Thus the secondary branch-cut is the green contour inside the blue region.  And in Figure \ref{figure:figure35b} the cut domain completely encloses the cut trace for branch $\{-1,m\}$.  And the cut trace slices through successive negative leaf sheets (once per sheet) in increasingly larger spiral arcs. The following is a function for the branch-cut that is used to  construct a derived $\pLog$ stack:
\begin{equation}
\text{branchCutF}(\theta,n)=e^{a/b(\theta+2 n\pi)}
\end{equation}
where as before $a=Log[|z|]$ and $b=\Arg[z]$. 

In the case of $n>0$, the secondary branch-cut is the intersection of
\begin{equation}
\begin{aligned}
r&>\text{exp}\bigg\{\frac{-\Arg(z)}{\ln(z)}\left(\Arg(w)+2 n\pi)\right)\bigg\}, \\
r&=\text{exp}\bigg\{\frac{\ln(z)}{\Arg(z)}\left(\Arg(w)+2n\pi)\right)\bigg\},
\end{aligned}
\end{equation}
and for the purpose of analyzing the branch-cuts, the quantity $(\Arg(w)+2n\pi)$ can be considered a constant. But for $z$ in Region 3A and $n>0$, $\displaystyle\frac{\ln(z)}{\Arg(z)}<\frac{-\Arg(z)}{\ln(z)}$ and therefore, the cut trace is always outside the cut domain so that the positive leaf sheets have only the single primary branch-cut.

%
%
\begin{figure}
     \centering
     \begin{subfigure}{0.32\textwidth}
         \centering
         \includegraphics[width=0.75\textwidth]{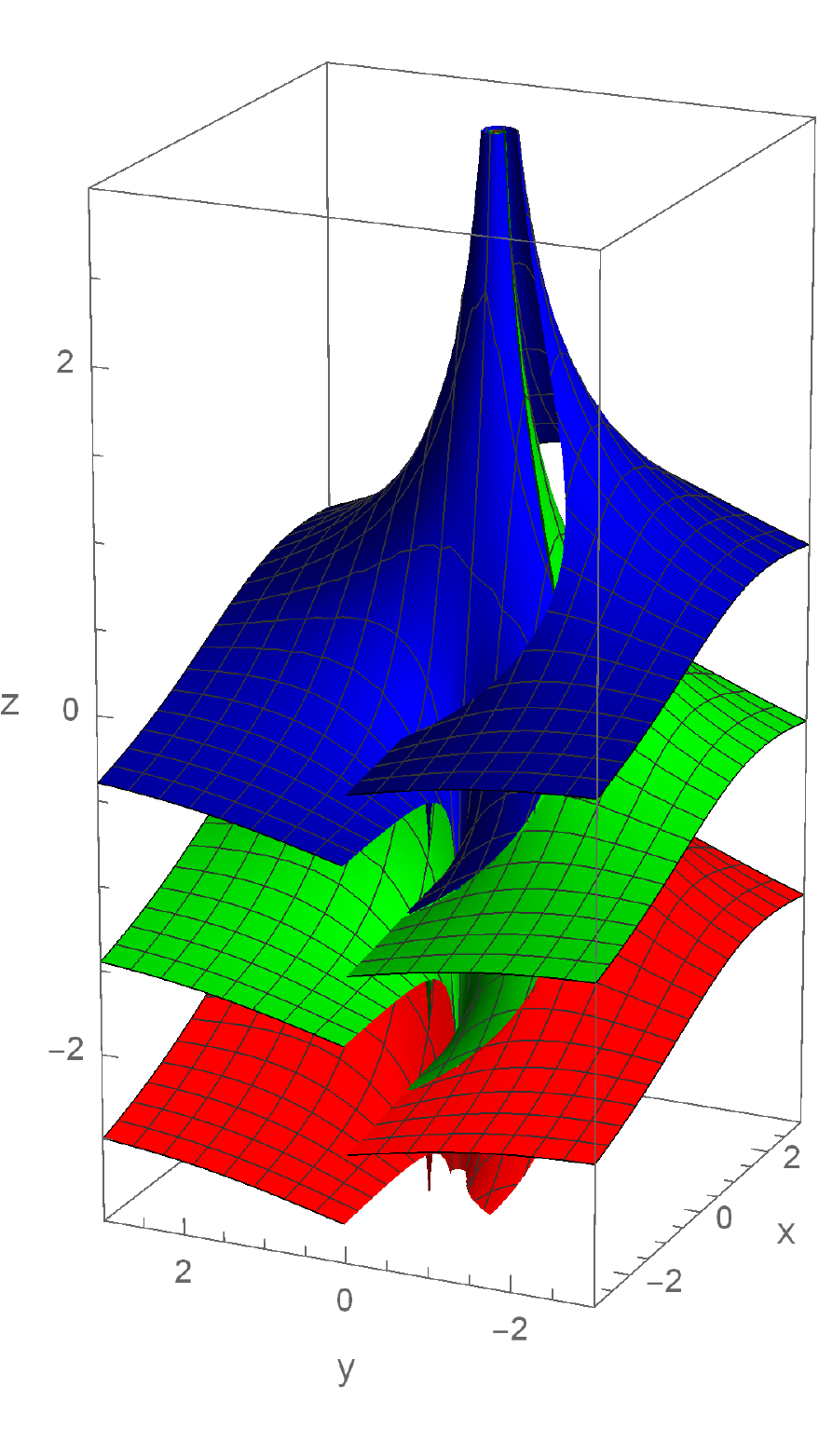}
         \caption{$\pLog$}
         \label{figure:figure36a}
     \end{subfigure}
     \hfill
		 \begin{subfigure}{0.32\textwidth}
         \centering
								\includegraphics[width=0.75\textwidth]{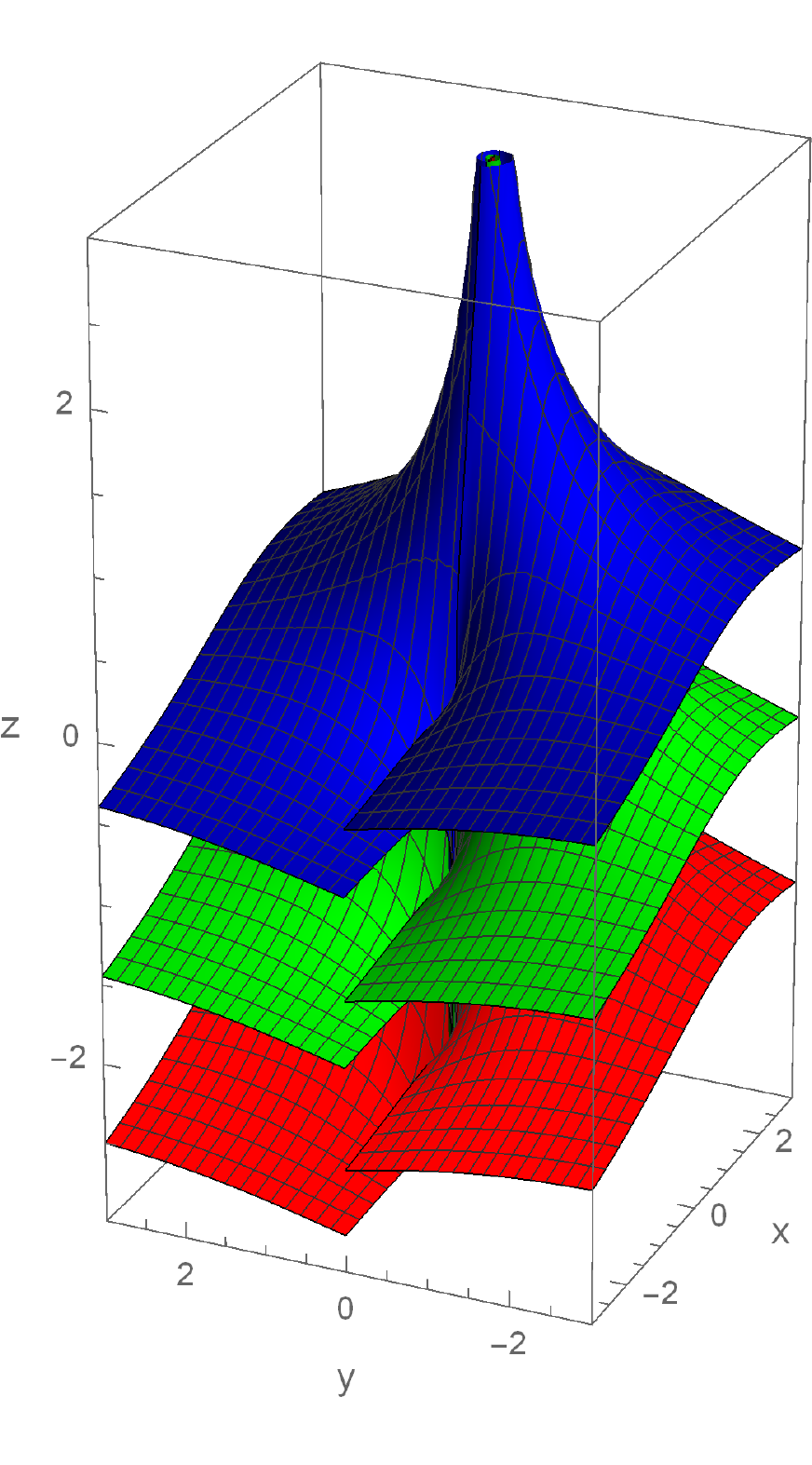}
         \caption{$\texttt{pLog3AN}$}
         \label{figure:figure36b}
     \end{subfigure}
		\hfill
     \begin{subfigure}{0.32\textwidth}
         \centering
         \includegraphics[width=0.75\textwidth]{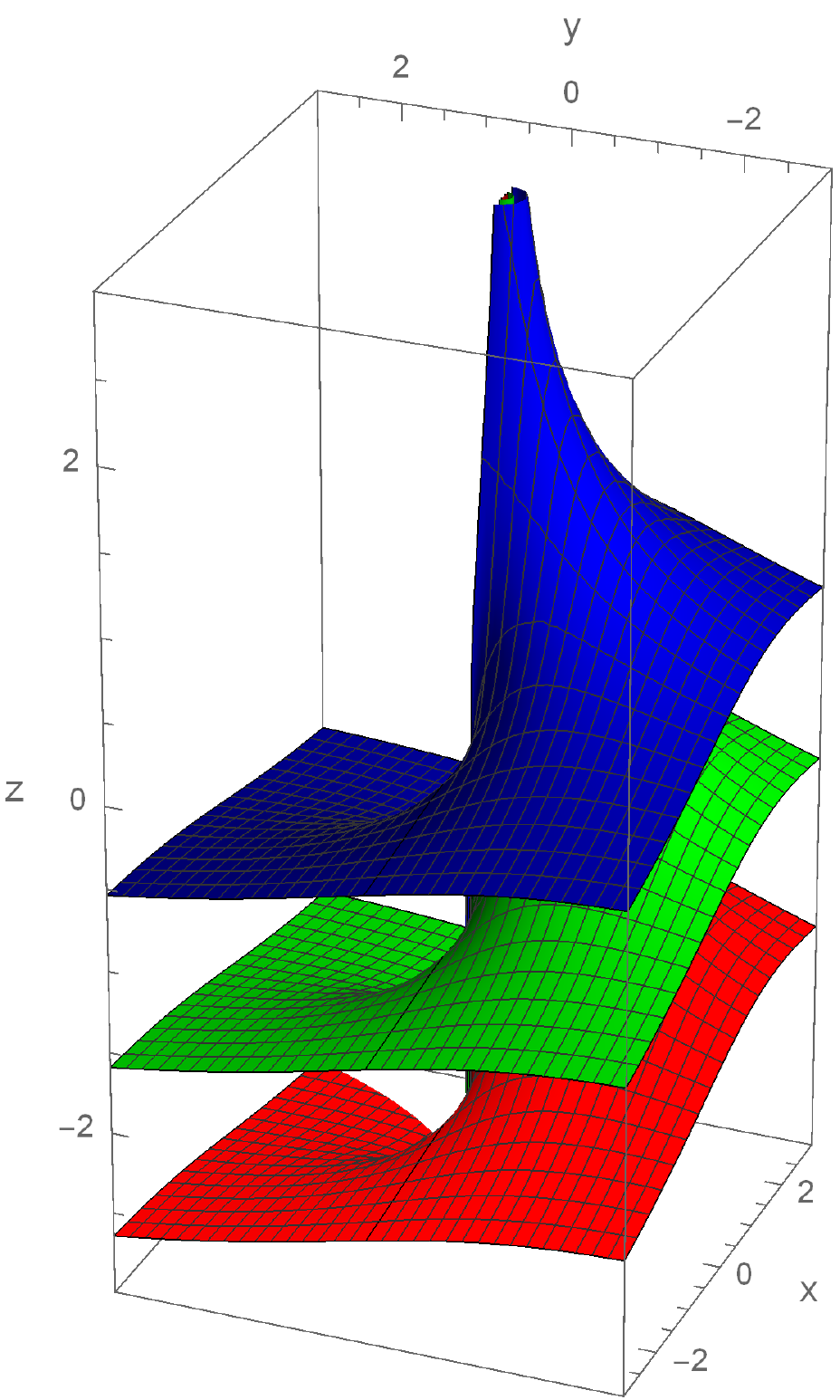}
         \caption{$\texttt{pLog3AP}$}
         \label{figure:figure36c}
     \end{subfigure}
     \hfill
     \caption{Region 3A default and derived log sheets $\{0,0\},\{-1,0\},\{1,0\}$}
        \label{figure:figure36}
\end{figure}

Figure \ref{figure:figure36} illustrates two types of log stacks derived from the default stack (the view perspective is looking down at the stacks from a point above and near the negative real axis).   Figure \ref{figure:figure36a} illustrates how the secondary branch-cut  passes through leaf sheets $\{0,m\}$ in the default case.  The small spiral arc of the branch-cut can be seen cutting through the sheets in the lower half-plane.  

Next, construct derived $\pLog$ sheets similar to what was done for Region 2A:  For the $\{0,m\}$ branches, the sheets are stitched with a single branch-cut $(-\infty,1)$.  For the negative-indexed sheets $\{-n,m\}$, the stitching produces a branch-cut $(-\infty,0)$.  Figure \ref{figure:figure36b} is a derived \texttt{pLog3AN} stack having a primary branch-cut along the negative real axis.  Figure \ref{figure:figure36c} is \texttt{pLog3AP} with the a single branch-cut along the positive real axis.  Listing \ref{code3} is the code to construct the \texttt{pLog3AP} stack.

\begin{center}
\begin{minipage}{0.75\linewidth}
\begin{lstlisting}[language=Mathematica,
frame=single,
caption=Mathematica code for the pLog3AP stack,
label=code3]
pLog3AP[w_, n_, m_] := Module[{rVal, aVal, pVal, cutVal, logVal},
   rVal = Abs[w];
   aVal = Arg[w];
   If[n >= 0,
    If[aVal >= 0,
      logVal = pLog[w, n, m];
      ,
      logVal = pLog[w, n + 1, m];
      ];
    ,
    If[aVal < 0,
      cutVal = branchCutF[aVal, n + 1];
      If[rVal >= cutVal,
       logVal = pLog[w, n + 1, m];
       ,
       logVal = pLog[w, n + 1, m - 1];
       ];
      ,
      cutVal = branchCutF[aVal, n];
      If[rVal >= cutVal,
       logVal = pLog[w, n, m];
       ,
       logVal = pLog[w, n, m - 1];
       ];
      ];
    ];
   logVal
   ];
	\end{lstlisting}
	\end{minipage}
	\end{center}

These particular stacks were derived because the default $\pLog$ iterator fails across the primary branch-cut. 
\begin{figure}[!ht]
     \centering
     \begin{subfigure}{0.4\textwidth}
         \centering
         \includegraphics[width=0.75\textwidth]{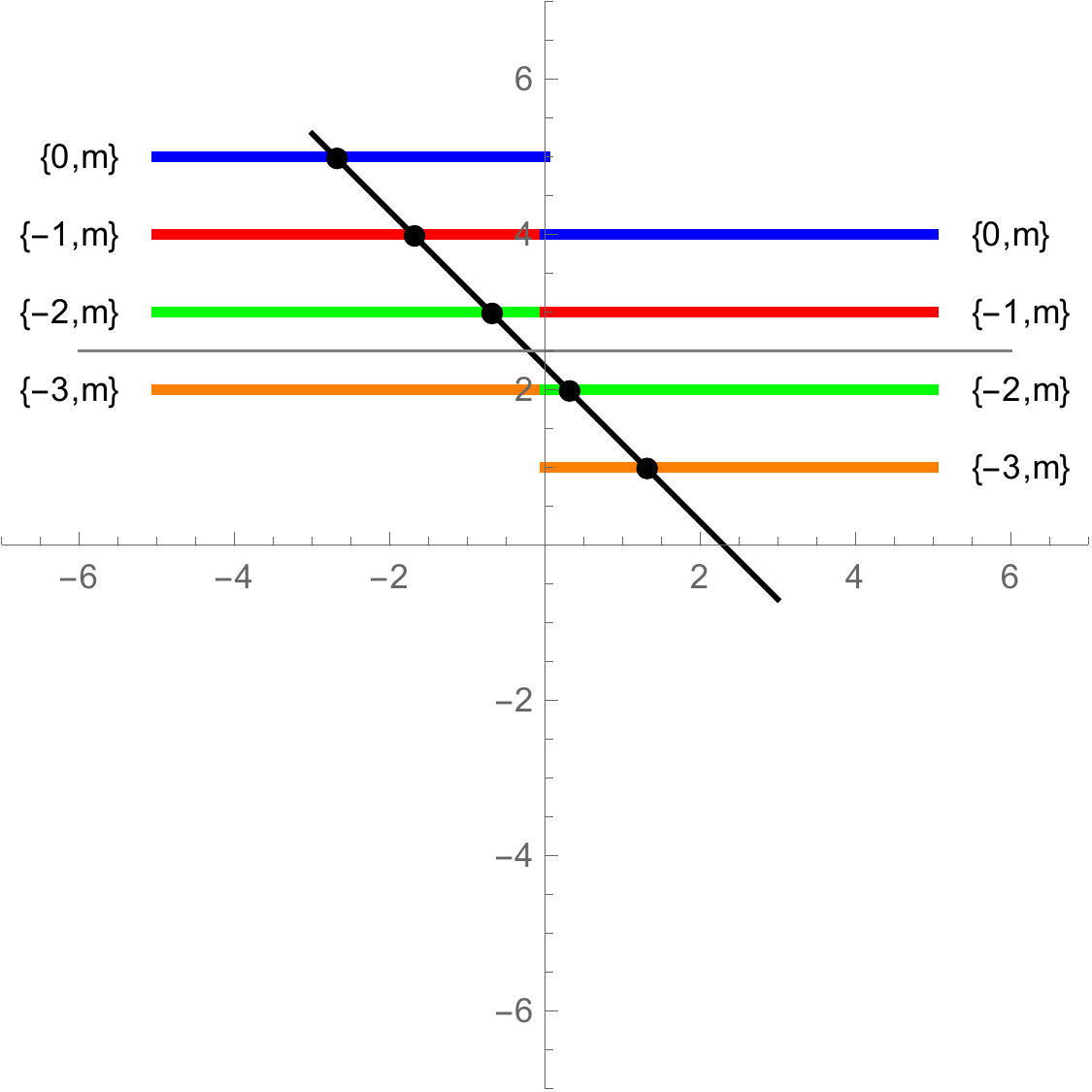}
         \caption{Negative leaf orientation $\{0,m\}$}
         \label{figure:figure351a}
     \end{subfigure}
     \hfill
     \begin{subfigure}{0.4\textwidth}
         \centering
         \includegraphics[width=0.75\textwidth]{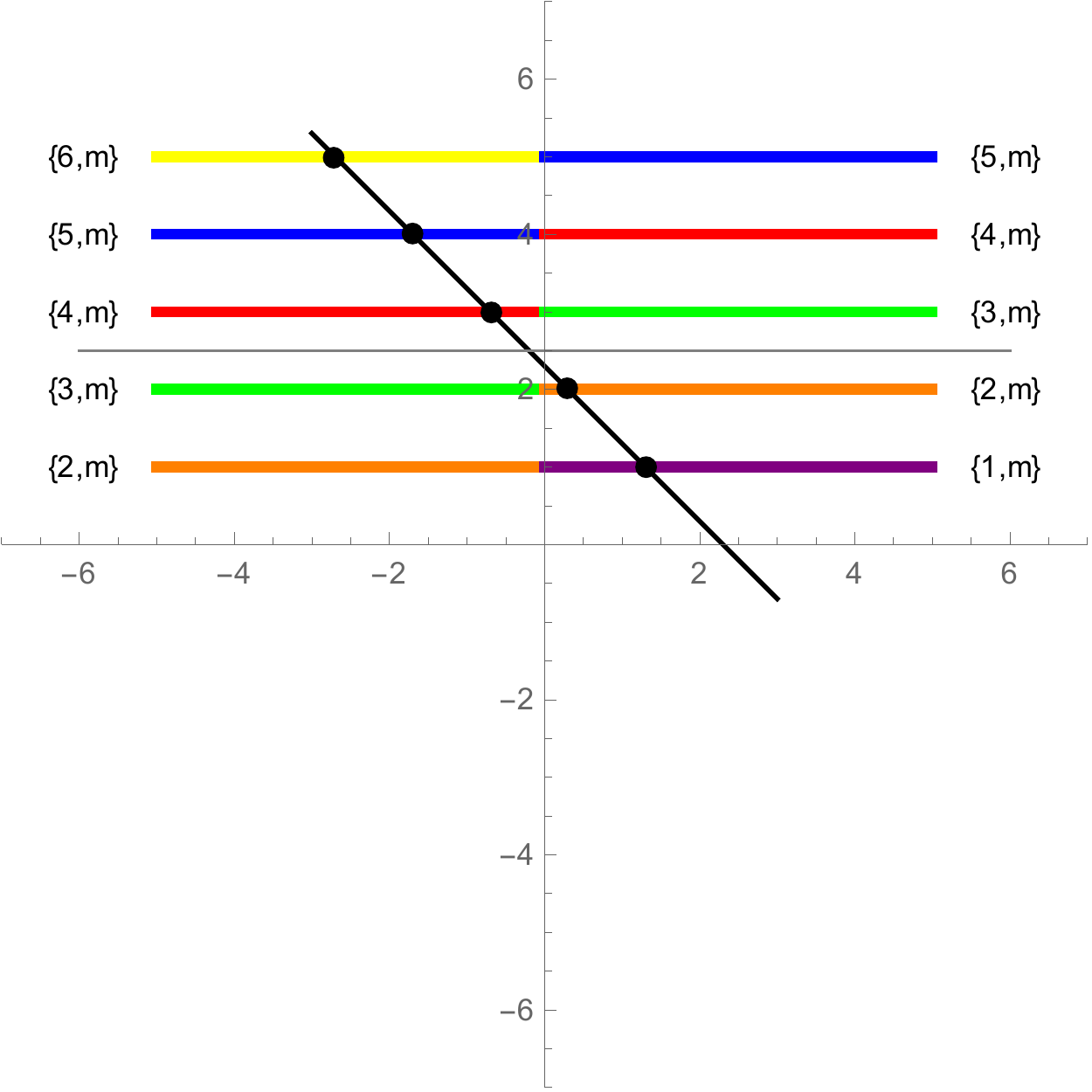}
         \caption{Positive leaf orientation $\{-1,m\}$}
         \label{figure:figure351b}
     \end{subfigure}
     \hfill
     \caption{Leaf sheets across primary branch-cut of $\pLog$}
        \label{figure:figure351}
\end{figure}

 Figure \ref{figure:figure351} shows a cross-section of the $\pLog$ leaf sheets over the primary branch-cut between different color sheets, and a black line representing the root trace across the sheets (where the line intersects a sheet is a root).  Notice in Figure \ref{figure:figure351a}  there are two roots on the green leaf.  An iteration of this sheet with $\pLog$ would miss one of the roots because of the branch-cut on the negative real axis.  And in Figure \ref{figure:figure351b} the root trace misses the green sheet entirely causing the $\pLog$ iterator to diverge for this $\{n,m\}$ index.   These effect can be eliminated by iterating the negative numbered branches with \texttt{pLog3AP} having the branch cut along the positive real axis.  Figure \ref{figure:figure36c} shows this stack. 

\begin{figure}[!ht]
	\centering
			\includegraphics[scale=.65]{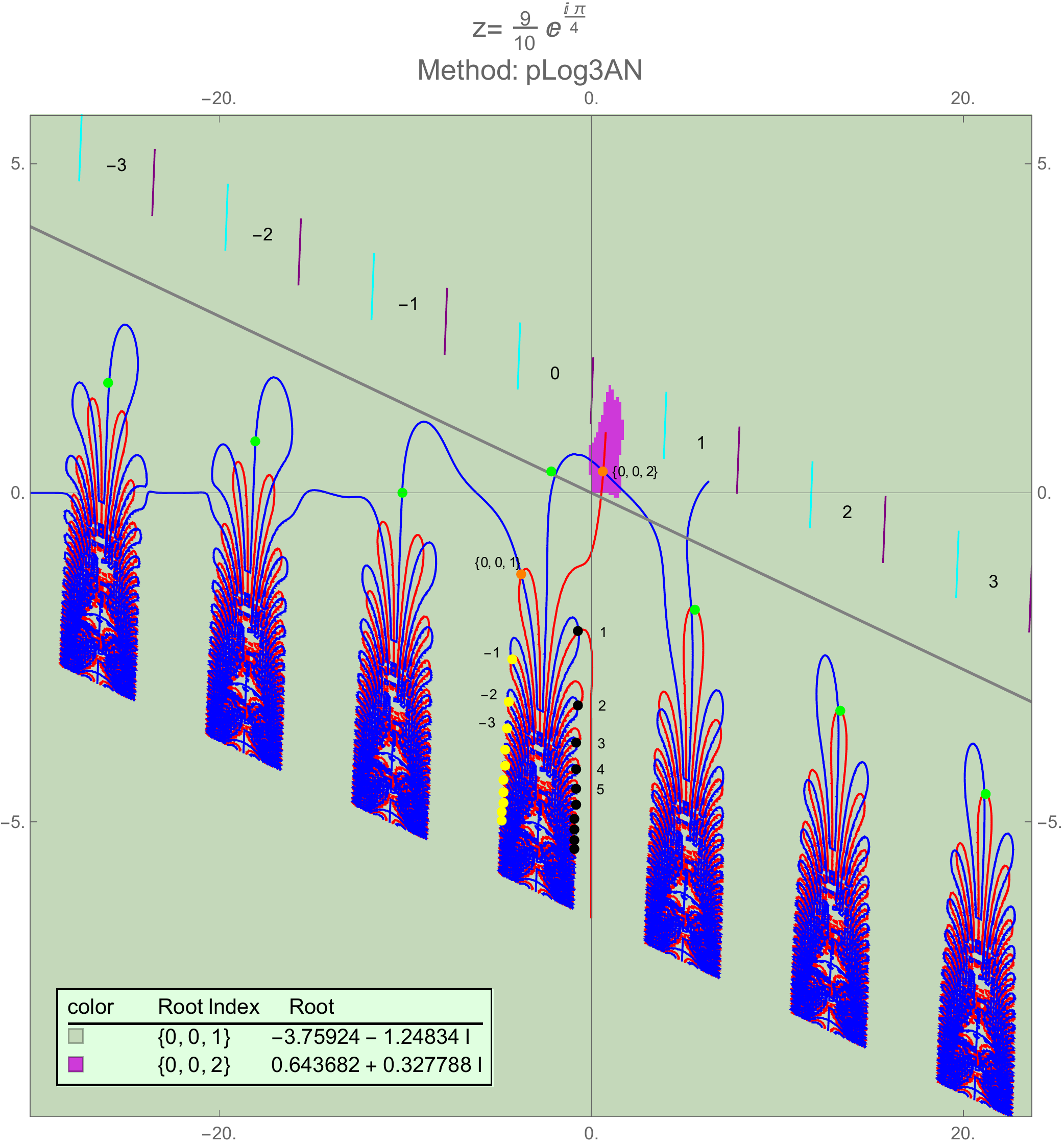}
		\caption{Region 3A contour plot with basins for root $\{0,0\}$}
   \label{figure:figure37}
\end{figure}
Figure \ref{figure:figure37} shows a \texttt{pLog3AN} basin diagram for root $\{0,0\}$.
The orange points are roots $\{0,0,1\}$ and $\{0,0,2\}$ for the two basins.   Note the small purple basin for $\{0,0,2\}$.  The roots depicted in the diagram were computed to an accuracy of $29$ digits in under $5$ iterations.
%
%
\section{Region 3B}
Region 3B is similar to Region3A except the inclination angle and branch-cuts are reversed. As $z$ approaches the origin, the angle becomes more positive.

The effect of the secondary branch-cuts of Region 3B is opposite that of Region 3A:  The cut trace slices through successive positive leaves in this region.  The negative traces spiral to the origin and has little effect on iterations.  And as with Region 3A, in order to prevent the effect of skipping roots when iterating across the principal branch-cut, both \texttt{pLog3BN} and \texttt{pLog3BP} are derived to compute the roots in this region. 
%
%
\section{Region 4A}
The branch inclination angle in this region is best analyzed on a half-circle of radius r.  Recall $\alpha=\arctan(a/b)$ where $\Log(z)=a+bi$. 
		\begin{figure}[!ht]
     \centering
     \begin{subfigure}{0.4\textwidth}
         \centering
         \includegraphics[width=0.75\textwidth]{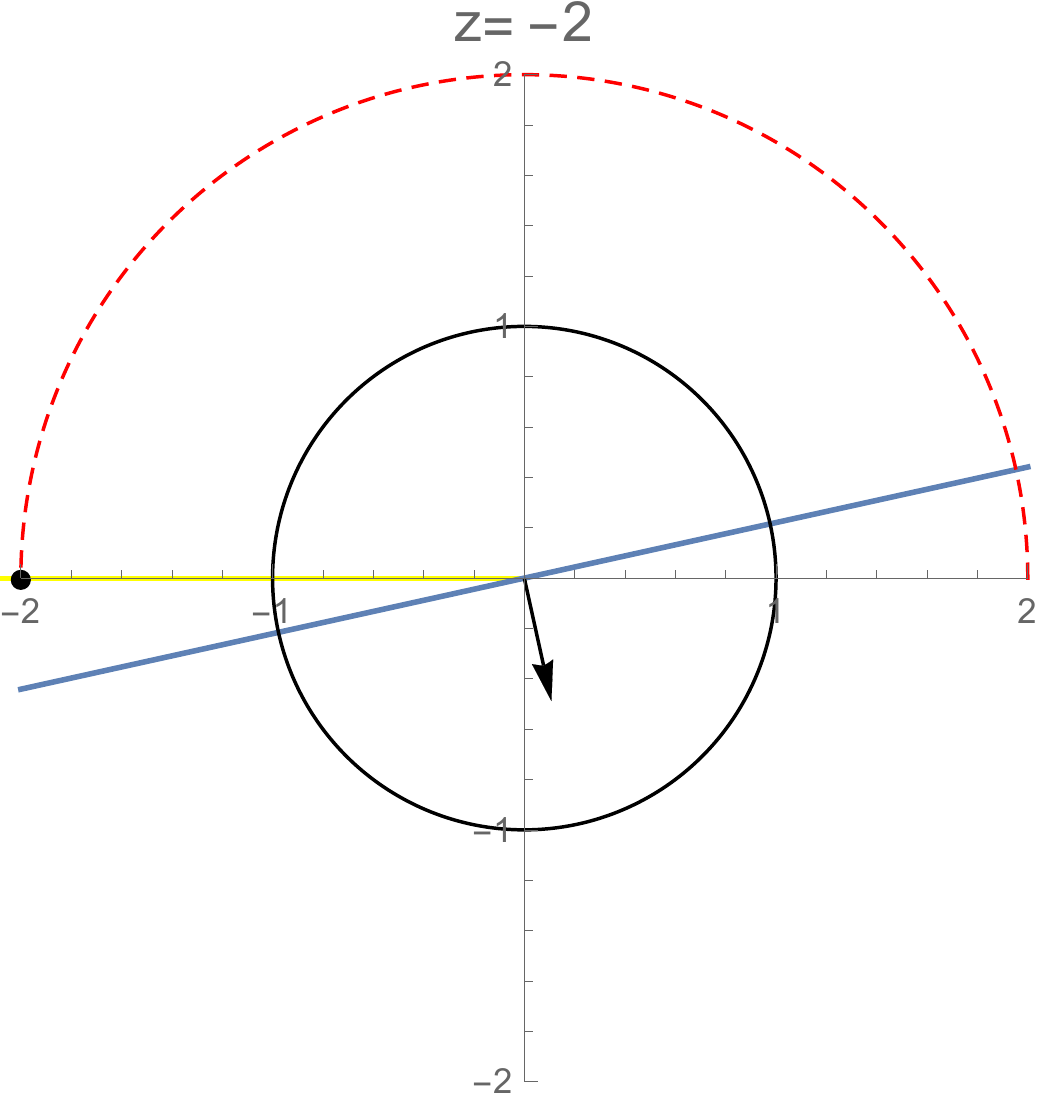}
         \caption{Inclination angle for $z=-2$}
         \label{figure:figure42a}
     \end{subfigure}
     \hfill
     \begin{subfigure}{0.4\textwidth}
         \centering
         \includegraphics[width=0.75\textwidth]{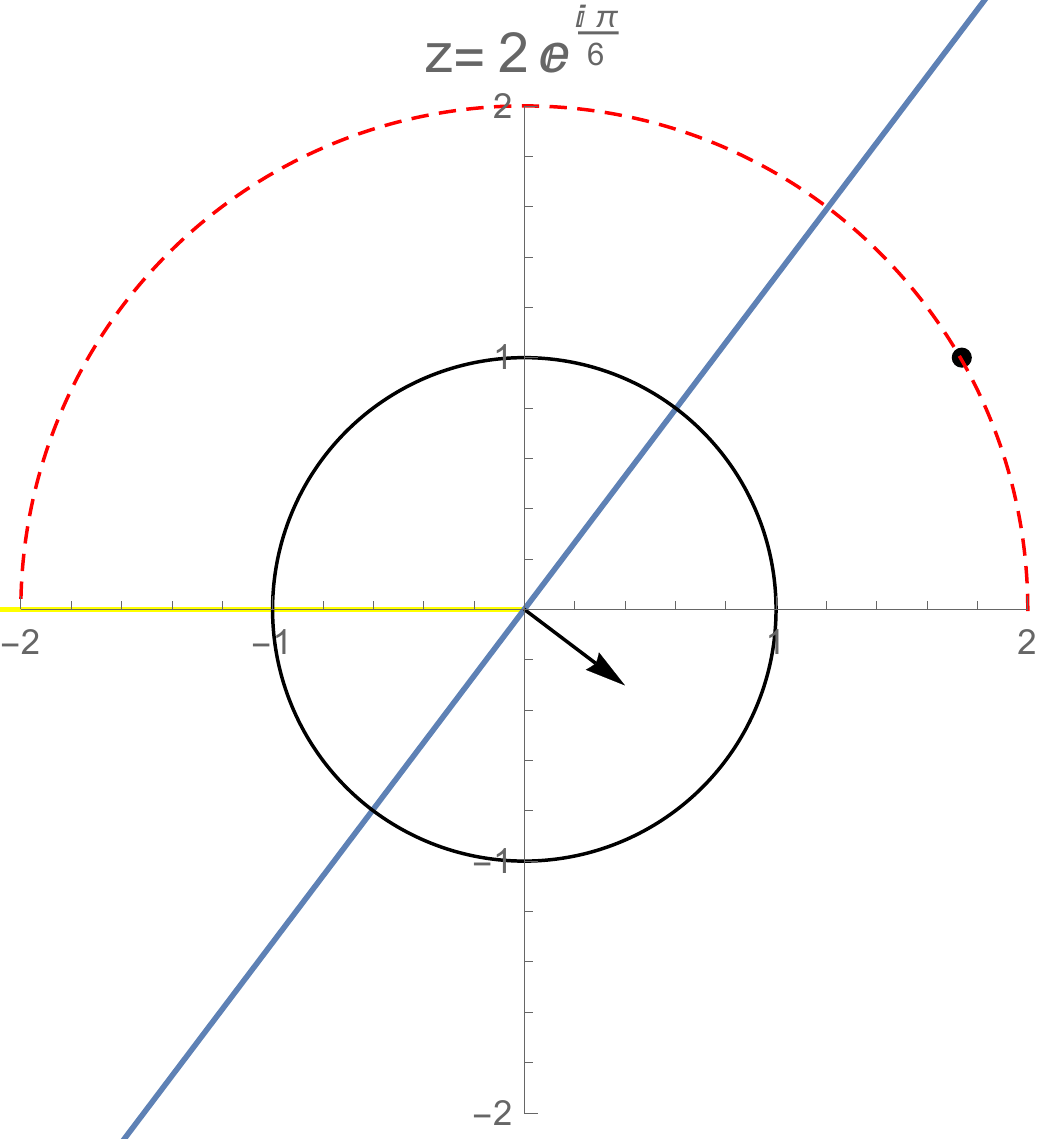}
         \caption{Inclination angle for $z=2e^{\pi i/6}$}
         \label{figure:figure42b}
     \end{subfigure}
     \hfill
     \caption{Example inclination angles for Region 4A}
        \label{figure:figure42}
\end{figure}

 Figure \ref{figure:figure42} illustrates how $\alpha$ changes as $z$ moves along the dashed red semi-circle.  Considering the expression $\alpha=\arctan(a/b)$, as $z$ moves clockwise along the half-circle, $a$ will remain constant while $b$ decreases.  This causes the argument of $\arctan$ to increase.  This effect becomes more pronounced as $|z|$ increases.   In this region, the contours open downward.

Only the non-negative leaf sheets are affected by the secondary branch-cut in this region.  This can be seen by analyzing the trace and cut functions:
\begin{equation*}
\begin{aligned}
a\ln(r)+bk&<0 \\
ak-b\ln(r)&=0
\end{aligned}
\end{equation*}  

Using the techniques described in previous sections, a derived \texttt{pLog4AN} stack with a single branch-cut along the negative real axis can be constructed to iterate the roots of this region since the branch line will eventually cross the positive real axis.

Table \ref{table:table4} shows the results for $z=10^{12}e^{\pi i/4}$ for the first $10$ roots over one-trillion on branch one-trillion.  After two iterations, the precision of the results was $58$ digits (the digits in Table \ref{table:table4} were truncated to conserve space).

\begin{table}[!ht]
\small
\centering
\begin{tabular}{|c|c|}
\hline
Root & Value \\
\hline
1&	$6.458402872023776511614937392*10^9+2.27212482053650764394928639099*10^{11}i$\\
2&	$6.458402872023776511614973554*10^9+2.27212482053650764394928638071*10^{11}i$\\
3&	$6.458402872023776511615009716*10^9+2.27212482053650764394928637043*10^{11}i$\\
4&	$6.458402872023776511615045878*10^9+2.27212482053650764394928636016*10^{11}i$\\
5&	$6.458402872023776511615082040*10^9+2.27212482053650764394928634988*10^{11}i$\\
6&	$6.458402872023776511615118202*10^9+2.27212482053650764394928633960*10^{11}i$\\
7&	$6.458402872023776511615154364*10^9+2.27212482053650764394928632932*10^{11}i$\\
8&	$6.458402872023776511615190526*10^9+2.27212482053650764394928631904*10^{11}i$\\
9&	$6.458402872023776511615226688*10^9+2.27212482053650764394928630876*10^{11}i$\\
10&	$6.458402872023776511615262850*10^9+2.27212482053650764394928629848*10^{11}i$\\
\hline
\end{tabular}
\caption{First $10$ roots above $10^{12}$ on branch $10^{12}$ for $z=10^{12}e^{\pi i/4}$}
 \label{table:table4}
\end{table}

\section{Region 4B}

\begin{figure}[!ht]
	\centering
			\includegraphics[scale=.65]{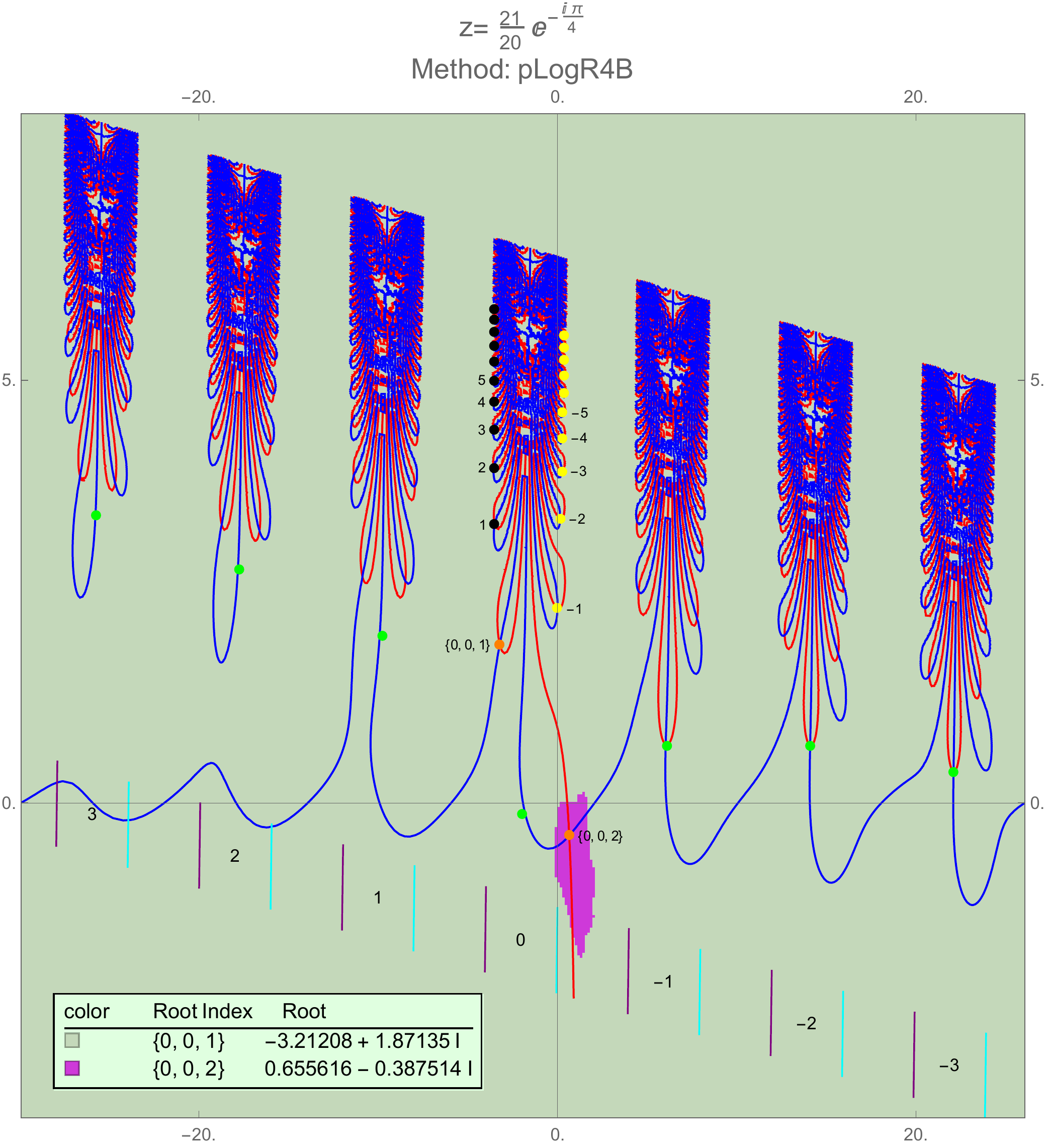}
		\caption{Contour and basin plot of region $4B$ example}
   \label{figure:figure44}
\end{figure}

Region 4B is the inverse of 4A. In this case, the secondary branch-cut affect the zero and positive leaf branches.  It should now be a simple matter to construct derived \texttt{pLog4BN} and \texttt{pLog4BP} stacks.  Figure \ref{figure:figure44} is an example contour and basin diagram showing two basins for root $\{0,0\}$.

\section{Conclusions}

The algorithms described in this paper provide a simple method to compute numerically-accessible roots of $T_2$, linking a root ID, associated contour plot, and the underlying logarithmic geometry of the function.  Most of the test cases studied dealt with branches near the origin  because the zeros near the origin were most challenging to compute, often producing multiple basins.  None of the branches tested away from the origin exhibited multiple basins.  Test cases were chosen to  best represent each of the $12$ regions defined above, or best illustrate the methods.  Roots as high as $\{10^{12},10^{12}\}$ were computed relatively easily and with high precision.  However it is likely roots extremely close to critical points such as the transitions between regions, or roots with extremely high values of $n$ and $m$ would pose greater challenges requiring further adjustments to the $\pLog$ iterator.  

The author's web-site, \href{http://www.jujusdiaries.com}{Algebraic Functions} will describe in detail the software used in this paper and further work done on the subject.

%
%

\end{document}